\documentclass[seceqn,MSNbibl,number,citesort,dvips]{arxbj}
\usepackage{mathbh}

\aid{0}
\volume{20}
\issue{4}
\pubyear{2014}
\firstpage{2247}
\lastpage{2277}
\doi{10.3150/13-BEJ556} %kopijuoti is PTS

\makeatletter

\newcommand{\rright}{\right}
\newcommand{\lleft}{\left}
\newcommand{\rrVert}{\Vert}
\newcommand{\rrvert}{\vert}
\newcommand{\llVert}{\Vert}
\newcommand{\llvert}{\vert}
\newcommand{\AR}{\operatorname{AR}}

\newcommand{\CC}{\mathsf{C}}
\newcommand{\DD}{\mathsf{D}}
\newcommand{\NN}{\mathbb{N}}
\newcommand{\RR}{\mathbb{R}}
\newcommand{\ZZ}{\mathbb{Z}}

\newcommand{\bM}{{\mathbf{M}}}
\newcommand{\bP}{{\mathbf{P}}}
\newcommand{\bQ}{{\mathbf{Q}}}
\newcommand{\bR}{{\mathbf{R}}}
\newcommand{\tbu}{\widetilde{{\mathbf{u}}}}
\newcommand{\bV}{{\mathbf{V}}}
\newcommand{\obV}{{\overline{\bV}}}
\newcommand{\bX}{{\mathbf{X}}}
\newcommand{\bZ}{{\mathbf{Z}}}
\newcommand{\bU}{{\mathbf{U}}}
\newcommand{\bgamma}{{\boldsymbol{\gamma}}}
\newcommand{\bxi}{{\boldsymbol{\xi}}}
\newcommand{\bzeta}{{\boldsymbol{\zeta}}}
\newcommand{\bvare}{{\boldsymbol{\vare}}}
\newcommand{\bzero}{{\mathbf{0}}}
\newcommand{\bone}{{\mathbf{1}}}

\newcommand{\cA}{{\mathcal A}}
\newcommand{\cB}{{\mathcal B}}
\newcommand{\cD}{{\mathcal D}}
\newcommand{\cF}{{\mathcal F}}
\newcommand{\cM}{{\mathcal M}}
\newcommand{\bcM}{\mathbf{\cM}}
\newcommand{\cN}{{\mathcal N}}
\newcommand{\bcN}{\mathbf{\cN}}
\newcommand{\cP}{{\mathcal P}}
\newcommand{\bcP}{\mathbf{\cP}}
\newcommand{\cU}{{\mathcal U}}
\newcommand{\bcU}{\mathbf{\cU}}
\newcommand{\cX}{{\mathcal X}}
\newcommand{\bcX}{\mathbf{\cX}}
\newcommand{\cY}{{\mathcal Y}}
\newcommand{\cZ}{{\mathcal Z}}
\newcommand{\cW}{{\mathcal W}}
\newcommand{\bcW}{\mathbf{\cW}}
\newcommand{\tbcW}{\widetilde{\bcW}}
\newcommand{\bcY}{\mathbf{\cY}}
\newcommand{\bcZ}{\mathbf{\cZ}}
\newcommand{\tcW}{\widetilde{\cW}}
\newcommand{\ttcW}{\widetilde{\tcW}}

\newcommand{\slu}{{\mathrm{lu}}}

\newcommand{\EE}{\mathrm{\mathbb{E}}}
\newcommand{\PP}{\mathrm{\mathbb{P}}}
\newcommand{\OO}{\mathrm{O}}
\newcommand{\var}{\operatorname{Var}}

\newcommand{\halpha}{\widehat{\alpha}}
\newcommand{\hbeta}{\widehat{\beta}}
\newcommand{\hdelta}{\widehat{\delta}}
\newcommand{\hvarrho}{\widehat{\varrho}}
\newcommand{\tH}{\widetilde{H}}
\newcommand{\vare}{\varepsilon}

\newcommand{\stoch}{\stackrel{\PP}{\longrightarrow}}
\newcommand{\distr}{\stackrel{\cD}{\longrightarrow}}
\newcommand{\distre}{\stackrel{\cD}{=}}
\newcommand{\lu}{\stackrel{\slu}{\longrightarrow}}
\newcommand{\as}{\stackrel{{\mathrm{a.s.}}}{\longrightarrow}}
\newcommand{\ase}{\stackrel{{\mathrm{a.s.}}}{=}}

\newcommand{\bbone}{\mathbh{1}}

\newtheorem{Thm}{Theorem}[section]
\newtheorem{Lem}[Thm]{Lemma}
\newtheorem{Pro}[Thm]{Proposition}
\newtheorem{Cor}[Thm]{Corollary}

\newremark{Rem}[Thm]{Remark}
\newcommand{\eqref}[1]{(\ref{#1})}
\newcommand{\binom}[2]{\pmatrix{#1\cr#2}}
\makeatother

\begin{document}
\begin{frontmatter}

\title{Asymptotic behavior of CLS estimators for $2$-type doubly symmetric
critical Galton--Watson processes with immigration}

\runtitle{$2$-type critical Galton--Watson processes with immigration}

\begin{aug}
%%%% inicialai - be tarpu
\author[1]{\inits{M.}\fnms{M\'arton} \snm{Isp\'any}\thanksref{1}\ead[label=e1]{ispany.marton@inf.unideb.hu}},
\author[2]{\inits{K.}\fnms{Krist\'of} \snm{K\"ormendi}\corref{}\thanksref{2,e2}\ead[label=e2,mark]{kormendi@math.u-szeged.hu}}
\and
\author[2]{\inits{G.}\fnms{Gyula} \snm{Pap}\thanksref{2,e3}\ead[label=e3,mark]{papgy@math.u-szeged.hu}}
%%\runauthor{} %% auto
\address[1]{University of Debrecen, Faculty of Informatics, Department of Information Technology, Pf. 12, H-4010 Debrecen, Hungary.
\printead{e1}}
\address[2]{University of Szeged, Faculty of Science, Bolyai Institute, Department of Stochastics, Aradi v\'ertan\'uk tere~1, H-6720 Szeged, Hungary.
\printead{e2};\\
\printead*{e3}}
\end{aug}

% HISTORY:
\received{\smonth{12} \syear{2012}}
\revised{\smonth{8} \syear{2013}}

% ABSTRACT
%
\begin{abstract}
In this paper, the asymptotic behavior of the conditional least squares (CLS)
estimators of the offspring means $(\alpha, \beta)$ and of the criticality
parameter $\varrho:= \alpha+ \beta$ for a $2$-type critical doubly
symmetric positively regular Galton--Watson branching process with immigration
is described.
\end{abstract}

% KEYWORDS
% visi is mazosios raides ir pagal abecele
%

\end{frontmatter}

%s1 #&#
\section{Introduction}
\label{section_intro}

Asymptotic behavior of CLS estimators for critical Galton--Watson
processes is
available only for single-type processes, see Wei and Winnicki \cite{WW,WW2} and Winnicki \cite{Win}, see also the monograph of Guttorp
\cite{Gut}.
In the present paper, the asymptotic behavior of the CLS estimators of the
offspring means and criticality parameter for $2$-type critical doubly
symmetric positively regular Galton--Watson process with immigration is
described, see Theorem~\ref{main}.
This study can be considered as the first step of examining the asymptotic
behavior of the CLS estimators of parameters of multitype critical branching
processes with immigration.
Shete and Sriram \cite{SheSri} obtained convergence results for
weighted CLS
estimators in the supercritical case.

Let us recall the results for a single-type Galton--Watson branching process
$(X_k)_{k \in\ZZ_+}$ with immigration and with initial value $X_0 = 0$.
Suppose that it is critical, that is, the offspring mean equals 1.
Wei and Winnicki \cite{WW} proved a functional limit theorem
$\cX^{(n)} \distr\cX$ as $n \to\infty$, where
$\cX^{(n)}_t := n^{-1} X_{{\lfloor nt\rfloor}}$ for $t \in\RR_+$,
$n \in\NN$, where
$\lfloor x \rfloor$ denotes the (lower) integer part of $x \in\RR$,
and $(\cX_t)_{t \in\RR_+}$ is a (nonnegative) diffusion process with
initial value $\cX_0 = 0$ and with generator
\[
Lf(x) = m_\vare f'(x) + \tfrac{1}{2}
V_\xi x f''(x) ,\qquad  f \in
C^\infty_{\mathrm{c}}(\RR_+) ,
\]
where $m_\vare$ denotes the immigration mean, $V_\xi$ denotes the
offspring variance, and $C^\infty_{\mathrm{c}}(\RR_+)$ denotes the
space of
infinitely differentiable functions on $\RR_+$ with compact support.
The process $(\cX_t)_{t\in\RR_+}$ can also be characterized as the unique
strong solution of the stochastic differential equation (SDE)
\[
\mathrm{d}\cX_t = m_\vare\,\mathrm{d}t + \sqrt{
V_\xi\cX_t^+ } \,\mathrm{d}\cW_t ,\qquad  t \in\RR_+ ,
\]
with initial value $\cX_0 = 0$, where $(\cW_t)_{t\in\RR_+}$ is a
standard Wiener process, and $x^+$ denotes the positive part of
$x \in\RR$.
Note that this so-called square-root process is also known as Feller
diffusion, or Cox--Ingersoll--Ross model in financial mathematics
(see Musiela and Rutkowski \cite{MR}, page 290).
In fact, $(4 V_\xi^{-1} \cX_t)_{t\in\RR_+}$ is the square of a
$4 V_\xi^{-1} m_\vare$-dimensional Bessel process started at 0 (see
Revuz and
Yor \cite{RevYor}, XI.1.1).

Assuming that the immigration mean $m_\vare$ is known, for the conditional
least squares estimator (CLSE)
\[
\halpha_n(X_1, \dots, X_n) =
\frac{\sum_{k=1}^n X_{k-1} (X_k - m_\vare)}{\sum_{k=1}^n X_{k-1}^2}
\]
of the offspring mean based on the observations $X_1, \dots, X_n$, one can
derive
\[
n\bigl(\halpha_n(X_1, \dots, X_n) - 1\bigr)
\distr \frac{\int_0^1 \cX_t \,\mathrm{d}(\cX_t - m_\vare t)}{\int_0^1 \cX
_t^2 \,\mathrm{d}t} \qquad \mbox{as }n \to\infty.
\]
(Wei and Winnicki \cite{WW2} contains a similar result for the CLS
estimator of
the offspring mean when the immigration mean is unknown.)

In Section~\ref{section_preliminaries}, we recall some preliminaries on
$2$-type Galton--Watson models with immigration.
Section~\ref{section_main_results} contains our main results.
Sections~\ref{section_proof}, \ref{section_proof_main},
\ref{section_proof_main1} and \ref{section_proof_maint} contain the proofs.
Appendix~\ref{section_estimators} is devoted to the CLS estimators.
In Appendix~\ref{section_moments}, we present estimates for the
moments of the
processes involved.
Appendices \ref{app_C} and \ref{section_conv_step_processes} are for
a version of
the continuous mapping theorem and for convergence of random step processes,
respectively.
For a detailed discussion of the whole paper, see Isp\'any \textit{et al.} \cite{IKP}.

%s2 #&#
\section{Preliminaries on $2$-type Galton--Watson models
with immigration}
\label{section_preliminaries}

Let $\ZZ_+$, $\NN$, $\RR$ and $\RR_+$ denote the set of
nonnegative integers, positive integers, real numbers and non-negative real
numbers, respectively.
Every random variable will be defined on a fixed probability space
$(\Omega,\cA,\PP)$.

For each $k, j \in\ZZ_+$ and $i, \ell\in\{ 1, 2 \}$, the number of
individuals of type $i$ in the $k$th generation will be
denoted by $X_{k,i}$, the number of type $\ell$ offsprings produced by
the $j$th individual who is of type $i$ belonging to the
$(k-1)$th generation will be denoted by $\xi_{k,j,i,\ell
}$, and
the number of type $i$ immigrants in the $k$th generation
will be denoted by $\vare_{k,i}$.
Then
%
%e2.1 #&#
\begin{equation}
\label{GWI(2)} %
\lleft[\matrix{ X_{k,1}
\cr
X_{k,2} }
\rright] %
= \sum_{j=1}^{X_{k-1,1}}
\lleft[\matrix{ \xi_{k,j,1,1}
\cr
\xi_{k,j,1,2} } \rright]
+\sum_{j=1}^{X_{k-1,2}} %
\lleft[\matrix{ \xi_{k,j,2,1}
\cr
\xi_{k,j,2,2} } \rright] %
+ %
\lleft[\matrix{ \vare_{k,1}
\cr
\vare_{k,2} }
\rright] %
,\qquad  k \in\NN.
\end{equation}
Here
$ \{ \bX_0, \bxi_{k,j,i}, \bvare_k
: k, j \in\NN, i \in\{ 1, 2 \}  \}$
are supposed to be independent, where
\[
\bX_k := %
\lleft[\matrix{ X_{k,1}
\cr
X_{k,2} } \rright] %
,\qquad  \bxi_{k,j,i} := %
\lleft[\matrix{ \xi_{k,j,i,1}
\cr
\xi_{k,j,i,2} } \rright] %
, \qquad \bvare_k := %
\lleft[\matrix{ \vare_{k,1}
\cr
\vare_{k,2} } \rright] %
.
\]
Moreover, $\{ \bxi_{k,j,1} : k, j \in\NN\}$,
$\{ \bxi_{k,j,2} : k, j \in\NN\}$ and $\{ \bvare_k : k \in\NN\}$
are supposed to consist of identically distributed random vectors.

We suppose $\EE(\|\bxi_{1,1,1}\|^2) < \infty$,
$\EE(\|\bxi_{1,1,2}\|^2) < \infty$ and $\EE(\|\bvare_1\|^2) <
\infty$.
Introduce the notations
\begin{eqnarray*}
{\mathbf{m}}_{\bxi_i} &:=& \EE (\bxi_{1,1,i} ) \in\RR
^2_+ ,\qquad  {\mathbf{m}}_{\bxi} := %
\lleft[\matrix{ {
\mathbf{m}}_{\bxi_1}  {\mathbf{m}}_{\bxi_2} } \rright] %
\in\RR^{2 \times2}_+ ,
\\
\bV_{\bxi_i} &:=& \var (\bxi_{1,1,i} ) \in\RR^{2 \times
2} ,\qquad
\obV_\bxi:= \tfrac{1}{2} (\bV_{\bxi_1} +
\bV_{\bxi_2}) \in\RR ^{2 \times2} ,
\\
{\mathbf{m}}_{\bvare} &:=& \EE (\bvare_1 ) \in
\RR^2_+ ,\qquad  \bV_{\bvare} := \var (\bvare_1 ) \in
\RR^{2 \times2} .
\end{eqnarray*}
Note that many authors define the offspring mean matrix as
${\mathbf{m}}^\top_\bxi$.
For $k \in\ZZ_+$, let
$\cF_k := \sigma ( \bX_0, \bX_1 , \dots, \bX_k  )$.
By \eqref{GWI(2)},
%
%e2.2 #&#
\begin{equation}
\label{mart} \EE(\bX_k |\cF_{k-1}) = X_{k-1,1} {
\mathbf{m}}_{\bxi_1} + X_{k-1,2} {\mathbf {m}}_{\bxi_2}
+{\mathbf{m}}_{\bvare} = {\mathbf{m}}_{\bxi} \bX_{k-1} +
{\mathbf{m}}_{\bvare} .
\end{equation}
Consequently, $\EE(\bX_k) = {\mathbf{m}}_{\bxi} \EE(\bX
_{k-1}) + {\mathbf{m}}_{\bvare}$,
$k \in\NN$, which implies
%
%e2.3 #&#
\begin{equation}
\label{EXk} \EE(\bX_k) = {\mathbf{m}}_{\bxi}^k
\EE(\bX_0) + \sum_{j=0}^{k-1} {
\mathbf{m}}_{\bxi}^j {\mathbf {m}}_{\bvare} ,\qquad  k
\in\NN.
\end{equation}
Hence, the offspring mean matrix ${\mathbf{m}}_{\bxi}$ plays a
crucial role in the
asymptotic behavior of the sequence $(\bX_k)_{ k \in\ZZ_+ }$.
Since ${\mathbf{m}}_{\bxi}$ has nonnegative entries, the
Frobenius--Perron theorem
(see, e.g., Horn and Johnson \cite{HJ}, Theorems 8.2.11 and 8.5.1) describes
the behavior of the powers ${\mathbf{m}}_{\bxi}^k$ as $k \to
\infty$.
According to this behavior, a $2$-type Galton--Watson process
$(\bX_k)_{k \in\ZZ_+}$ with immigration is referred to respectively as
\emph{subcritical}, \emph{critical} or \emph{supercritical} if
$\varrho< 1$, $\varrho= 1$ or $\varrho> 1$, where $\varrho$
denotes the spectral radius of the offspring mean matrix ${\mathbf
{m}}_\bxi$
(see, e.g., Athreya and Ney \cite{AN} or Quine \cite{Q}).
We will consider doubly symmetric $2$-type Galton--Watson processes with
immigration, when the offspring mean matrix has the form
%
%e2.4 #&#
\begin{equation}
\label{bA} {\mathbf{m}}_\bxi:= %
\lleft[\matrix{ \alpha&
\beta
\cr
\beta& \alpha } \rright] %
.
\end{equation}
Its spectral radius is $\varrho= \alpha+ \beta$, which will be called
\emph{criticality parameter}.
We will focus only on \emph{positively regular} doubly symmetric $2$-type
Galton--Watson processes with immigration, that is, when there is a positive
integer $k \in\NN$ such that the entries of ${\mathbf{m}}_\bxi
^k$ are positive
(see Kesten and Stigum \cite{KesSti1}), which is equivalent with
$\alpha> 0$ and $\beta> 0$.

For the sake of simplicity, we consider a zero start Galton--Watson process
with immigration, that is, we suppose $\bX_0 = \bzero$.
In the sequel, we always assume ${\mathbf{m}}_\vare\ne\bzero$, otherwise
$\bX_k = \bzero$ for all $k \in\NN$.

%s3 #&#
\section{Main results}
\label{section_main_results}

In order to find CLS estimators of the criticality parameter
$\varrho= \alpha+ \beta$, we introduce a further parameter
$\delta:= \alpha- \beta$.
Then $\alpha= (\varrho+ \delta)/2$ and $\beta= (\varrho- \delta)/2$,
thus the recursion \eqref{regr} can be written in the form
\[
\bX_k = \frac{1}{2} %
\lleft[\matrix{ \varrho+
\delta& \varrho- \delta
\cr
\varrho- \delta& \varrho+ \delta } \rright] %
\bX_{k-1} + \bM_k + {\mathbf{m}}_\bvare, \qquad k \in\NN.
\]
For each $n \in\NN$, a CLS estimator $(\hspace*{0.5pt}\hvarrho_n, \hdelta_n)$ of
$(\varrho, \delta)$ based on a sample $\bX_1, \ldots, \bX_n$ can be
obtained by minimizing the sum of squares
\[
\sum_{k=1}^n \biggl\llVert
\bX_k - \frac{1}{2} %
\left[\matrix{ \varrho+ \delta&
\varrho- \delta
\cr
\varrho- \delta& \varrho+ \delta }\right] %
\bX_{k-1} - {
\mathbf{m}}_\bvare\biggr\rrVert ^2
\]
with respect to $(\varrho, \delta)$ over $\RR^2$, and it has the form
%
%e3.1 #&#
%e3.2 #&#
\begin{eqnarray}
\label{CLSEr}\hvarrho_n &:=& \frac{\sum_{k=1}^n
\langle\bone, {\mathbf{x}}_k - {\mathbf{m}}_\bvare\rangle
\langle\bone, {\mathbf{x}}_{k-1} \rangle} {
\sum_{k=1}^n \langle\bone, {\mathbf{x}}_{k-1} \rangle^2} ,
\\
\label{CLSEd}\hdelta_n &:=& \frac{\sum_{k=1}^n
\langle\tbu, {\mathbf{x}}_k - {\mathbf{m}}_\bvare\rangle
\langle\tbu, {\mathbf{x}}_{k-1} \rangle} {
\sum_{k=1}^n \langle\tbu, {\mathbf{x}}_{k-1} \rangle^2}
\end{eqnarray}
on the set $H_n \cap\tH_n$, where
\[
\bone:= %
\left[\matrix{ 1
\cr
1 }\right] %
\in\RR^2 , \qquad \tbu:= %
\left[\matrix{ 1
\cr
-1 }\right] %
\in\RR^2 ,
\]
and
%
%e3.3 #&#
%e3.4 #&#
\begin{eqnarray}
\label{H_n}H_n &:=& \Biggl\{ ({\mathbf{x}}_1, \ldots, {
\mathbf{x}}_n) \in\bigl(\RR^2\bigr)^n : \sum
_{k=1}^n \langle\bone, {
\mathbf{x}}_{k-1} \rangle^2 > 0 \Biggr\} ,
\\
\label{tH_n}\tH_n &:=& \Biggl\{ ({\mathbf{x}}_1, \ldots, {
\mathbf{x}}_n) \in\bigl(\RR^2\bigr)^n : \sum
_{k=1}^n \langle\tbu, {
\mathbf{x}}_{k-1} \rangle^2 > 0 \Biggr\} ,
\end{eqnarray}
where ${\mathbf{x}}_0 := \bzero$ is the zero vector in $\RR^2$.
In a natural way, we extend the CLS estimators $\hvarrho_n$ and
$\hdelta_n$ to the set $H_n$ and $\tH_n$, respectively.
Moreover, for each $n \in\NN$, any CLS estimator $(\halpha_n, \hbeta_n)$
of the offspring means $(\alpha, \beta)$ based on a sample
$\bX_1, \ldots, \bX_n$ has the form
%
%e3.5 #&#
\begin{equation}
\label{CLSEab} %
\left[\matrix{ \halpha_n
\cr
\hbeta_n }\right] %
= \frac{1}{2} %
\left[\matrix{ 1 & 1
\cr
1 & -1 }\right] %
\left[\matrix{ \hvarrho_n
\vspace*{2pt}\cr
\hdelta_n }\right] %
,
\end{equation}
whenever the sample belongs to the set $H_n \cap\tH_n$.
For the proof see Isp\'any \textit{et al.} \cite{IKP}, Lemma A.1.

In what follows, we always assume that $(\bX_k)_{k \in\ZZ_+}$ is a
$2$-type doubly symmetric Galton--Watson process with offspring means
$(\alpha, \beta) \in(0, 1)^2$ such that $\alpha+ \beta= 1$
(hence it is critical and positively regular), $\bX_0 = \bzero$,
$\EE(\|\bxi_{1,1,1}\|^8) < \infty$, $\EE(\|\bxi_{1,1,2}\|^8) <
\infty$,
$\EE(\|\bvare_1\|^8) < \infty$, and ${\mathbf{m}}_\bvare\ne
\bzero$.
Then $\lim_{n \to\infty} \PP( (\bX_1, \ldots, \bX_n) \in H_n) = 1$.
If $\langle\obV_\bxi\tbu, \tbu\rangle> 0$, or if
$\langle\obV_\bxi\tbu, \tbu\rangle= 0$ and
$\EE(\langle\tbu, \bvare_1 \rangle^2) > 0$, then
$\lim_{n \to\infty} \PP( (\bX_1, \ldots, \bX_n) \in\tH_n) =
1$, see
Proposition~\ref{ExUn}.

Let $(\cY_t)_{t \in\RR_+}$ be the unique strong solution of the stochastic
differential equation (SDE)
%
%e3.6 #&#
\begin{equation}
\label{Y} \mathrm{d}\cY_t = \langle\bone, {\mathbf{m}}_\bvare
\rangle\,\mathrm{d}t + \sqrt{\langle\obV_\bxi\bone, \bone\rangle
\cY_t^+} \,\mathrm {d}\cW_t , \qquad t \in\RR_+ ,\
\cY_0 = 0 ,
\end{equation}
where $(\cW_t)_{t \in\RR_+}$ is a standard Wiener process.
%
%th3.1 #&#
\begin{Thm}\label{main}
We have
%
%e3.7 #&#
\begin{equation}
\label{rho} n (\hspace*{0.5pt}\hvarrho_n -1 ) \distr \frac{\int_0^1 \cY_t \,\mathrm{d}(\cY_t - \langle\bone,
{\mathbf{m}}_\bvare\rangle t)} {
\int_0^1 \cY_t^2 \,\mathrm{d}t} ,\qquad  \mbox{as }n \to\infty.
\end{equation}

If $\langle\obV_\bxi\bone, \bone\rangle= 0$, then
%
%e3.8 #&#
\begin{equation}
\label{rho1} n^{3/2} (\hspace*{0.5pt}\hvarrho_n -1 ) \distr \cN \biggl(0,
\frac{3 \langle\bV_\bvare\bone, \bone\rangle} {
\langle\bone, {\mathbf{m}}_\bvare\rangle^2} \biggr) , \qquad \mbox{ as }n \to\infty.
\end{equation}

If $\langle\obV_\bxi\tbu, \tbu\rangle> 0$, then
%
%e3.9 #&#
\begin{equation}
\label{alpha,beta} %
\left[\matrix{ n^{1/2} (\halpha_n -
\alpha)
\cr
n^{1/2} (\hbeta_n - \beta) } \right]%
\distr \sqrt{
\alpha\beta} \frac{\int_0^1 \cY_t \,\mathrm{d}\tcW_t} {
\int_0^1 \cY_t \,\mathrm{d}t} %
\left[\matrix{ 1
\cr
-1 }\right] %
, \qquad \mbox{as }n \to\infty,
\end{equation}
where $(\tcW_t)_{t \in\RR_+}$ is a standard Wiener process,
independent from
$(\cW_t)_{t \in\RR_+}$.

If $\langle\obV_\bxi\tbu, \tbu\rangle= 0$ and
$\EE(\langle\tbu, \bvare_1 \rangle^2) > 0$, then
%
%e3.10 #&#
\begin{equation}
\label{alpha,betat} %
\lleft[\matrix{ n^{1/2} (
\halpha_n - \alpha)
\cr
n^{1/2} (\hbeta_n -
\beta) } \rright] %
\distr \cN \biggl(0, \frac{\langle\bV_\bvare\tbu, \tbu\rangle} {
4 \EE(\langle\tbu, \bvare_1 \rangle^2)} \biggr)
\lleft[\matrix{ 1
\cr
-1 } \rright] %
, \qquad \mbox{as }n \to
\infty.
\end{equation}
\end{Thm}
%
%re3.2 #&#
\begin{Rem}\label{REMARK-1}
If $\langle\obV_\bxi\tbu, \tbu\rangle> 0$ and
$\langle\obV_\bxi\bone, \bone\rangle= 0$ then in \eqref{alpha,beta}
we have
\[
\sqrt{\alpha\beta} \frac{\int_0^1 \cY_t \,\mathrm{d}\tcW_t} {
\int_0^1 \cY_t \,\mathrm{d}t} %
\lleft[\matrix{ 1
\cr
-1 }
\rright] %
\distre \cN \biggl(0, \frac{4}{3}\alpha\beta \biggr)
\lleft[\matrix{ 1
\cr
-1 } \rright] %
.
\]
\end{Rem}
%
%re3.3 #&#
\begin{Rem}\label{REMARK0}
Note that the assumption $\langle\obV_\bxi\bone, \bone\rangle= 0$ is
fulfilled if and only if $\xi_{1,1,1,1} + \xi_{1,1,1,2} \ase1$ and
$\xi_{1,1,2,1} + \xi_{1,1,2,2} \ase1$, that is, the total number of
offsprings
produced by an individual of type 1 is 1, and the same holds for individuals
of type 2.
In a similar way, the assumption
$\langle\obV_\bxi\tbu, \tbu\rangle= 0$ is fulfilled if and only if
$\alpha= \beta= \frac{1}{2}$, $\xi_{1,1,1,1} \ase\xi_{1,1,1,2}$ and
$\xi_{1,1,2,1} \ase\xi_{1,1,2,2}$, that is, the number of offsprings
of type 1
and of type 2 produced by an individual of type 1 are the same, and the same
holds for individuals of type 2.
Observe that the assumptions $\langle\obV_\bxi\bone, \bone\rangle
= 0$
and $\langle\obV_\bxi\tbu, \tbu\rangle= 0$ can not be fulfilled at
the same time.

Condition $\EE(\langle\tbu, \bvare_1 \rangle^2) > 0$ fails to hold
if and
only if $\vare_{1,1} - \vare_{1,2} \ase0$, and, under the assumption
$\langle\obV_\bxi\tbu, \tbu\rangle= 0$, this implies
$X_{k,1} \ase X_{k,2}$ (see Lemma~\ref{main_VVt}), when
$\PP((\bX_1, \ldots, \bX_n) \in H_n \cap\tH_n) = 0$ for all
$n \in\NN$, and hence the LSE of the offspring means $(\alpha, \beta)$
is not defined uniquely, see Appendix~\ref{section_estimators}.
\end{Rem}

%re3.4 #&#
\begin{Rem}\label{REMARK1}
For each $n \in\NN$, consider the random step process
\[
\bcX^{(n)}_t := n^{-1} \bX_{\lfloor nt\rfloor},\qquad  t \in
\RR_+ .
\]
Theorem~\ref{main_conv} implies convergence \eqref{convXZ}, hence
%
%e3.11 #&#
\begin{equation}
\label{convX} \bcX^{(n)} \distr\bcX:= \tfrac{1}{2} \cY\bone\qquad  \mbox{as }n \to\infty,
\end{equation}
where the process $(\cY_t)_{t \in\RR_+}$ is the unique strong
solution of
the SDE \eqref{Y} with initial value $\cY_0 = 0$.
Note that convergence \eqref{convX} holds even if
$\langle\obV_\bxi\bone, \bone\rangle= 0$, when the unique strong
solution of \eqref{Y} is the deterministic function $\cY_t = \langle
\bone,
{\mathbf{m}}_\bvare\rangle t$, $t\in\RR_+$.

The SDE \eqref{Y} has a unique strong solution $(\cY_t^{(y)})_{t \in
\RR_+}$
for all initial values $\cY_0^{(y)} = y \in\RR$, and if $y \geq0$,
then $\cY_t^{(y)}$ is nonnegative for all $t \in\RR_+$ with
probability one, hence $\cY_t^+$ may be replaced by $\cY_t$ under the
square root in \eqref{Y}, see, for example,
Barczy \textit{et al.} \cite{BarIspPap0}, Remark~3.3.
\end{Rem}
%
%re3.5 #&#
\begin{Rem}\label{REMARK3}
We note that in the critical positively regular case the limit distributions
for the CLS estimators of the offspring means $(\alpha, \beta)$ are
concentrated on the line $\{(u, v) \in\RR^2 : u + v = 0\}$.
In order to handle the difficulty caused by this degeneracy, we use an
appropriate reparametrization.
Surprisingly, the scaling factor of the CLS estimators of $(\alpha,
\beta)$
is always $\sqrt{n}$, which is the same as in the subcritical case.
The reason of this strange phenomenon can be understood from the joint
asymptotic behavior of the numerator and the denominator of the CLS estimators
given in Theorems \ref{main_Ad}, \ref{main1_Ad} and \ref{maint_Ad}.
The scaling factor of the estimators of the criticality parameter
$\varrho$
is usually $n$, except in a particular special case of
$\langle\obV_\bxi\bone, \bone\rangle= 0$, when it is $n^{3/2}$.
One of the decisive tools in deriving the needed asymptotic behavior is a
good bound for the moments of the involved processes, see Corollary~\ref{EEX_EEU_EEV}.
\end{Rem}
%
%re3.6 #&#
\begin{Rem}\label{REMARK4}
The shape of
$\int_0^1 \cY_t \,\mathrm{d}(\cY_t - \langle\bone, {\mathbf
{m}}_\bvare\rangle t)
/ \int_0^1 \cY_t^2 \,\mathrm{d}t$
in \eqref{rho} is similar to the limit distribution of the Dickey--Fuller
statistics for unit root test of $\AR(1)$ time series, see, for exmple,
Hamilton \cite{Ham},\vadjust{\goodbreak} formulas 17.4.2 and 17.4.7, or Tanaka
\cite{Tan}, (7.14) and Theorem~9.5.1.
The shape of $\int_0^1 \cY_t \,\mathrm{d}\tcW_t  / \int_0^1 \cY
_t \,\mathrm{d}t$ in
\eqref{alpha,beta} is also similar, but it contains two independent
standard Wiener processes.
This phenomenon is very similar to the appearance of two independent standard
Wiener processes in limit theorems for CLS estimators of the variance
of the
offspring and immigration distributions for critical branching
processes with
immigration in Winnicki \cite{Win}, Theorems 3.5 and 3.8.
Finally, note that the limit distribution of the CLS estimator of the
criticality parameter $\varrho$ is non-symmetric and non-normal in
\eqref{rho}, and symmetric normal in \eqref{rho1}, but the limit distribution
of the CLS estimator of the offspring means $(\alpha, \beta)$ is always
symmetric, although non-normal in \eqref{alpha,beta}.
\end{Rem}
%
%re3.7 #&#
\begin{Rem}\label{REMARK5}
The eighth order moment conditions on the offspring and immigration
distributions in Theorem~\ref{main} seem
to be too strong, but we note that the process
$(\bX_k)_{k\in\ZZ_+}$ can be considered as a heteroscedastic time series.
Indeed, $\bX_k = {\mathbf{m}}_\bxi\bX_{k-1} + {\mathbf
{m}}_\bvare+ \bM_k$, see \eqref{regr},
and by \eqref{Mcond},
$\EE( \bM_k \bM_k^\top|\cF_{k-1} )
= X_{k-1,1} \bV_{\bxi_1} + X_{k-1,2} \bV_{\bxi_2} + \bV_\bvare$,
$k \in\NN$.
That is why we think that the behavior of the process $(\bX_k)_{k\in
\ZZ_+}$
is similar to GARCH models, where, even in the stable case, high moment
conditions are needed for convergence of estimators such as the quasi-maximum
likelihood estimator in Hall and Yao \cite{HalYao} or the Whittle estimator
in Mikosch and Straumann \cite{MikStr}.
\end{Rem}

%s4 #&#
\section{Proof of the main results}
\label{section_proof}

Applying \eqref{mart}, let us introduce the sequence
%
%e4.1 #&#
\begin{equation}
\label{Mk} \bM_k := \bX_k - \EE(\bX_k |
\cF_{k-1}) = \bX_k - {\mathbf{m}}_\bxi
\bX_{k-1} - {\mathbf{m}}_\bvare,\qquad  k \in\NN,
\end{equation}
of martingale differences with respect to the filtration
$(\cF_k)_{k \in\ZZ_+}$.
By \eqref{Mk}, the process $(\bX_k)_{k \in\ZZ_+}$ satisfies the recursion
%
%e4.2 #&#
\begin{equation}
\label{regr} \bX_k = {\mathbf{m}}_\bxi\bX_{k-1}
+ {\mathbf{m}}_\bvare+ \bM_k ,\qquad  k \in\NN.
\end{equation}
Next, let us introduce the sequence
\[
U_k := \langle\bone, \bX_k \rangle= X_{k,1} +
X_{k,2} ,\qquad  k \in\ZZ _+ .
\]
One can observe that $U_k \geq0$ for all $k \in\ZZ_+$, and
%
%e4.3 #&#
\begin{equation}
\label{rec_U} U_k = U_{k-1} + \langle\bone, {
\mathbf{m}}_\bvare\rangle + \langle\bone, \bM_k \rangle, \qquad k
\in\NN,
\end{equation}
since
$\langle\bone, {\mathbf{m}}_\bxi\bX_{k-1} \rangle= \bone
^\top{\mathbf{m}}_\bxi\bX_{k-1}
= \bone^\top\bX_{k-1} = U_{k-1}$,
because $\varrho= \alpha+ \beta= 1$ implies that $\bone$ is a
left eigenvector of the mean matrix ${\mathbf{m}}_\bxi$ belonging
to the eigenvalue
1.
Hence, $(U_k)_{k \in\ZZ_+}$ is a nonnegative unstable $\AR(1)$ process with
positive drift $\langle\bone, {\mathbf{m}}_\bvare\rangle$ and with
heteroscedastic innovation $(\langle\bone, \bM_k \rangle)_{k \in
\NN}$.
Moreover, let
\[
V_k := \langle\tbu, \bX_k \rangle= X_{k,1} -
X_{k,2} ,\qquad  k \in\ZZ_+ .
\]
Note that we have
%
%e4.4 #&#
\begin{equation}
\label{rec_V} V_k = (\alpha- \beta) V_{k-1} + \langle\tbu,
{\mathbf {m}}_\bvare\rangle + \langle\tbu, \bM_k
\rangle,\qquad  k \in\NN,\vadjust{\goodbreak}
\end{equation}
since
$\langle\tbu, {\mathbf{m}}_\bxi\bX_{k-1} \rangle= \tbu^\top
{\mathbf{m}}_\bxi\bX_{k-1}
= (\alpha- \beta) \tbu^\top\bX_{k-1} = (\alpha- \beta) V_{k-1}$,
because $\tbu$ is a left eigenvector of the mean matrix ${\mathbf
{m}}_\bxi$
belonging to the eigenvalue $\alpha- \beta$.
Thus $(V_k)_{k \in\NN}$ is a stable $\AR(1)$ process with drift
$\langle\tbu, {\mathbf{m}}_\bvare\rangle$ and with
heteroscedastic innovation
$(\langle\tbu, \bM_k \rangle)_{k \in\NN}$.
Observe that
%
%e4.5 #&#
\begin{equation}
\label{XUV} X_{k,1} = (U_k + V_k)/2 ,\qquad
X_{k,2} = (U_k - V_k)/2 ,\qquad  k \in\ZZ_+ .
\end{equation}

By \eqref{CLSEr}, for each $n \in\NN$, we have
\[
\hvarrho_n - 1 = \frac{\sum_{k=1}^n \langle\bone, \bM_k\rangle U_{k-1}} {
\sum_{k=1}^n U_{k-1}^2} ,
\]
whenever $(\bX_1, \ldots, \bX_n) \in H_n$, where $H_n$, $n \in\NN$,
are given in \eqref{H_n}.
By \eqref{CLSEd}, for each $n \in\NN$, we have
%
%e4.6 #&#
\begin{eqnarray}
\label{CLSE_delta-delta} \hdelta_n - \delta = \frac{\sum_{k=1}^n \langle\tbu, \bM_k\rangle V_{k-1}} {
\sum_{k=1}^n V_{k-1}^2} ,
\end{eqnarray}
whenever $(\bX_1, \ldots, \bX_n) \in\tH_n$, where $\tH_n$,
$n \in\NN$, are given in \eqref{tH_n}.

Theorem~\ref{main} will follow from the following statements by the continuous
mapping theorem.
%
%th4.1 #&#
\begin{Thm}\label{main_Ad}
We have, as $n \to\infty$,
\[
\sum_{k=1}^n %
\lleft[
\matrix{ n^{-3} U_{k-1}^2
\cr
n^{-2}
V_{k-1}^2
\cr
n^{-2} \langle\bone,
\bM_k \rangle U_{k-1}
\cr
n^{-3/2} \langle\tbu,
\bM_k \rangle V_{k-1} } \rright] %
\distr
\lleft[\matrix{ \displaystyle \int_0^1
\cY_t^2 \,\mathrm{d}t
\cr
(4 \alpha\beta)^{-1}
\langle\obV_\bxi\tbu, \tbu\rangle \int_0^1
\cY_t \,\mathrm{d}t
\cr
\displaystyle \int_0^1
\cY_t \,\mathrm{d}\bigl(\cY_t - \langle\bone, {
\mathbf{m}}_\bvare \rangle t\bigr)
\cr
(4 \alpha\beta)^{-1/2}
\langle\obV_\bxi\tbu, \tbu\rangle \displaystyle \int_0^1
\cY_t \,\mathrm{d}\tcW_t } \rright] %
.
\]
\end{Thm}
%
%th4.2 #&#
\begin{Thm}\label{main1_Ad}
If $\langle\obV_\bxi\bone, \bone\rangle= 0$ then, as $n \to
\infty$,
\[
\sum_{k=1}^n %
\lleft[
\matrix{ n^{-3} U_{k-1}^2
\cr
n^{-2}
V_{k-1}^2
\cr
n^{-3/2} \langle\bone,
\bM_k \rangle U_{k-1}
\cr
n^{-3/2} \langle\tbu,
\bM_k \rangle V_{k-1} } \rright] %
\distr
\lleft[\matrix{ \displaystyle \int_0^1
\cY_t^2 \,\mathrm{d}t
\cr
(4 \alpha\beta)^{-1}
\langle\obV_\bxi\tbu, \tbu\rangle \displaystyle \int_0^1
\cY_t \,\mathrm{d}t
\cr
\langle\bV_\bvare\bone, \bone
\rangle^{1/2} \int_0^1
\cY_t \,\mathrm{d}\ttcW_t
\cr
(4 \alpha\beta)^{-1/2}
\langle\obV_\bxi\tbu, \tbu\rangle \displaystyle \int_0^1
\cY_t \,\mathrm{d}\tcW_t } \rright] %
,
\]
where $(\ttcW_t)_{t\in\RR_+}$ is a standard Wiener process,
independent from
$(\cW_t)_{t \in\RR_+}$ and $(\tcW_t)_{t \in\RR_+}$.
Note that $(\cY_t)_{t \in\RR_+}$ is now the deterministic function
$\cY_t = \langle\bone, {\mathbf{m}}_\bvare\rangle t$, $t \in
\RR_+$, hence
$\int_0^1 \cY_t^2 \,\mathrm{d}t = \langle\bone, {\mathbf
{m}}_\bvare\rangle^2/3$,
$\int_0^1 \cY_t \,\mathrm{d}t = \langle\bone, {\mathbf
{m}}_\bvare\rangle/2$,
$\int_0^1 \cY_t \,\mathrm{d}\ttcW_t
= \langle\bone, {\mathbf{m}}_\bvare\rangle\int_0^1 t \,\mathrm
{d}\ttcW_t$
and
$\int_0^1 \cY_t \,\mathrm{d}\tcW_t
=\linebreak[4]  \langle\bone, {\mathbf{m}}_\bvare\rangle\int_0^1 t \,\mathrm
{d}\tcW_t$.\vspace*{-2pt}
\end{Thm}
%
%th4.3 #&#
\begin{Thm}\label{maint_Ad}
If $\langle\obV_\bxi\tbu, \tbu\rangle= 0$ then, as $n \to\infty$,
\begingroup
\abovedisplayskip=7pt
\belowdisplayskip=7pt
\[
\sum_{k=1}^n %
\lleft[
\matrix{ n^{-3} U_{k-1}^2
\cr
n^{-1}
V_{k-1}^2
\cr
n^{-2} \langle\bone,
\bM_k \rangle U_{k-1}
\cr
n^{-1/2} \langle\tbu,
\bM_k \rangle V_{k-1} } \rright] %
\distr
\lleft[\matrix{ \displaystyle \int_0^1
\cY_t^2 \,\mathrm{d}t
\cr
\EE\bigl(\langle\tbu,
\bvare_1 \rangle^2\bigr)
\cr
\displaystyle \int_0^1
\cY_t \,\mathrm{d}\bigl(\cY_t - \langle\bone, {
\mathbf{m}}_\bvare \rangle t\bigr)
\cr
\bigl[\langle\bV_\bvare
\tbu, \tbu\rangle \EE\bigl(\langle\tbu, \bvare_1 \rangle^2
\bigr)\bigr]^{1/2} \tcW_1 } \rright] %
.\vspace*{-2pt}
\]
\endgroup
\end{Thm}

%s5 #&#
\section{Proof of Theorem
\texorpdfstring{\protect\ref{main_Ad}}{4.1}}\vspace*{-2pt}
\label{section_proof_main}

Consider the sequence of stochastic processes
\begingroup
\abovedisplayskip=7pt
\belowdisplayskip=7pt
\[
\bcZ^{(n)}_t := %
\lleft[\matrix{
\bcM_t^{(n)}
\cr
\bcN_t^{(n)}
\cr
\bcP_t^{(n)} } \rright] %
:= \sum
_{k=1}^{\lfloor nt\rfloor} \bZ^{(n)}_k ,
\]
with
\[
\bZ^{(n)}_k := %
\lleft[\matrix{
n^{-1} \bM_k
\cr
n^{-2} \bM_k
U_{k-1}
\cr
n^{-3/2} \bM_k V_{k-1} }
\rright] %
= %
\lleft[\matrix{ n^{-1}
\cr
n^{-2} U_{k-1}
\cr
n^{-3/2} V_{k-1} } \rright] %
\otimes\bM_k
\]
\endgroup
for $t \in\RR_+$ and $k, n \in\NN$, where $\otimes$ denotes
Kronecker product of matrices.
Theorem~\ref{main_Ad} follows from Lemma~\ref{main_VV} and the following
theorem (this will be explained after Theorem~\ref{main_conv}).\vspace*{-2pt}
%
%th5.1 #&#
\begin{Thm}\label{main_conv}
We have
\begingroup
\abovedisplayskip=7pt
\belowdisplayskip=7pt
%e5.1 #&#
\begin{equation}
\label{conv_Z} \bcZ^{(n)} \distr\bcZ, \qquad \mbox{as }n\to\infty,
\end{equation}
where the process $(\bcZ_t)_{t \in\RR_+}$ with values in $(\RR
^2)^3$ is
the unique strong solution of the SDE
%
%e5.2 #&#
\begin{equation}
\label{ZSDE} \mathrm{d}\bcZ_t = \gamma(t, \bcZ_t)
\lleft[\matrix{ \mathrm{d}\bcW_t
\cr
\mathrm{d}
\tbcW_t } \rright] %
,\qquad  t \in\RR_+ ,
\end{equation}
with initial value $\bcZ_0 = \bzero$, where $(\bcW_t)_{t \in\RR
_+}$ and
$(\tbcW_t)_{t \in\RR_+}$ are independent $2$-dimensional standard Wiener
processes, and $\gamma\dvtx  \RR_+ \times(\RR^2)^3 \to(\RR^{2\times
2})^{3\times2}$
is defined by
\[
\gamma(t, {\mathbf{x}}) := %
\lleft[\matrix{ \bigl\langle\bone,
({\mathbf{x}}_1 + t {\mathbf{m}}_\bvare)^+ \bigr
\rangle^{1/2} \obV_\bxi^{1/2} & 0
\cr
\bigl\langle
\bone, ({\mathbf{x}}_1 + t {\mathbf{m}}_\bvare)^+ \bigr
\rangle^{3/2} \obV_\bxi^{1/2} & 0
\cr
0 & \biggl(
\displaystyle \frac{\langle\obV_\bxi\tbu, \tbu\rangle} {
4 \alpha\beta} \biggr)^{1/2} \langle\bone, {\mathbf{x}}_1
+ t {\mathbf{m}}_\bvare \rangle\obV_\bxi^{1/2} }
\rright] %
\]
\endgroup
for $t \in\RR_+$ and ${\mathbf{x}}= ({\mathbf{x}}_1 ,
{\mathbf{x}}_2 , {\mathbf{x}}_3) \in(\RR^2)^3$.\vadjust{\vspace*{-12pt}\goodbreak}
\end{Thm}
(Note that the statement of Theorem~\ref{main_conv} holds even if
$\langle\obV_\bxi\tbu, \tbu\rangle= 0$, when the last $2$-dimensional
coordinate process of the unique strong solution $(\bcZ_t)_{t\in\RR
_+}$ is
$\bzero$.)

The SDE \eqref{ZSDE} has the form
\begin{eqnarray*}
\mathrm{d}\bcZ_t = %
\lleft[\matrix{ \mathrm{d}
\bcM_t
\cr
\mathrm{d}\bcN_t
\cr
\mathrm{d}
\bcP_t } \rright] %
= %
\lleft[\matrix{
\bigl\langle\bone, (\bcM_t + t {\mathbf{m}}_\bvare)^+ \bigr
\rangle ^{1/2} \obV_\bxi^{1/2} \,\mathrm{d}
\bcW_t
\cr
\bigl\langle\bone, (\bcM_t + t {
\mathbf{m}}_\bvare)^+ \bigr\rangle ^{3/2} \obV_\bxi^{1/2}
\,\mathrm{d}\bcW_t
\cr
\biggl(\displaystyle \frac{\langle\obV_\bxi\tbu, \tbu\rangle} {
4 \alpha\beta} \biggr)^{1/2}
\langle\bone, \bcM_t + t {\mathbf{m}}_\bvare\rangle\obV
_\bxi^{1/2} \,\mathrm{d}\tbcW_t } \rright]
,\qquad  t\in\RR_+ .
\end{eqnarray*}
Isp\'any and Pap \cite{IspPap2} proved that the first $2$-dimensional equation
of this SDE has a unique strong solution $(\bcM_t)_{t\in\RR_+}$ with initial
value $\bcM_0 = \bzero$, and $(\bcM_t + t {\mathbf{m}}_\bvare
)^+$ may be
replaced by $\bcM_t + t {\mathbf{m}}_\bvare$ (see the proof of
\cite[Theorem~3.1]{IspPap2}).
Thus, the SDE \eqref{ZSDE} has a unique strong solution with initial value
$\bcZ_0 = \bzero$, and we have
\[
\bcZ_t = %
\lleft[\matrix{ \bcM_t
\cr
\bcN_t
\cr
\bcP_t } \rright] %
= %
\lleft[\matrix{ \displaystyle \int_0^t \langle\bone,
\bcM_t + t {\mathbf{m}}_\bvare\rangle^{1/2}
\obV_\bxi^{1/2} \,\mathrm{d}\bcW_s
\cr
\displaystyle \int
_0^t \langle\bone, \bcM_t + t {
\mathbf{m}}_\bvare\rangle\,\mathrm {d}\bcM_s
\cr
\biggl(
\displaystyle \frac{\langle\obV_\bxi\tbu, \tbu\rangle} {
4 \alpha\beta} \biggr)^{1/2} \int_0^t
\langle\bone, \bcM_t + t {\mathbf{m}}_\bvare\rangle\obV
_\bxi^{1/2} \,\mathrm{d}\tbcW_s } \rright]
,\qquad  t\in\RR_+ .
\]
By the method of the proof of $\cX^{(n)} \distr\cX$ in Theorem~3.1 in
Barczy \textit{et al.} \cite{BarIspPap0}, applying Lemma~\ref{Conv2Funct}, one can
easily derive
%
%e5.3 #&#
\begin{eqnarray}
\label{convXZ} %
\lleft[\matrix{ \bcX^{(n)}
\cr
\bcZ^{(n)} } \rright] %
\distr %
\lleft[\matrix{
\bcX
\cr
\bcZ } \rright] %
, \qquad \mbox{as }n \to\infty,
\end{eqnarray}
where
\[
\bcX^{(n)}_t := n^{-1} \bX_{\lfloor nt\rfloor},\qquad
\bcX_t: = \tfrac{1}{2} \langle\bone, \bcM_t + t {
\mathbf{m}}_\bvare\rangle\bone,\qquad  t \in\RR_+ ,\ n\in\NN,
\]
see Isp\'any \textit{et al.} \cite{IKP}, page 10.
Now, with the process
\[
\cY_t := \langle\bone, \bcX_t \rangle = \langle\bone,
\bcM_t + t {\mathbf{m}}_\bvare\rangle, \qquad t \in\RR_+ ,
\]
we have
\[
\bcX_t = \tfrac{1}{2} \cY_t \bone,\qquad  t \in\RR_+ .
\]
By It\^o's formula, we obtain that the process $(\cY_t)_{t\in\RR_+}$
satisfies
the SDE \eqref{Y}.
Next, similarly to the proof of \eqref{seged2}, by Lemma~\ref{Marci},
convergence \eqref{convXZ} and Lemma~\ref{main_VV} with
$U_{k-1} = \langle\bone, \bX_{k-1} \rangle$ implies
\[
\sum_{k=1}^n %
\lleft[
\matrix{ n^{-3} U_{k-1}^2
\cr
n^{-2}
V_{k-1}^2
\cr
n^{-2} \langle\bone,
\bM_k \rangle U_{k-1}
\cr
n^{-3/2} \langle\tbu,
\bM_k \rangle V_{k-1} } \rright] %
\distr
\lleft[\matrix{ \displaystyle \int_0^1
\langle\bone, \bcX_t \rangle^2 \,\mathrm{d}t
\cr
\displaystyle \frac{\langle\obV_\bxi\tbu, \tbu\rangle} {
4 \alpha\beta} \int_0^1 \langle\bone,
\bcX_t \rangle\,\mathrm{d}t
\cr
\displaystyle \int_0^1
\cY_t \,\mathrm{d}\langle\bone, \bcM_t \rangle
\cr
\biggl(
\displaystyle \frac{\langle\obV_\bxi\tbu, \tbu\rangle} {
4 \alpha\beta} \biggr)^{1/2} \int_0^1
\cY_t \,\mathrm{d}\bigl\langle\tbu, \obV_\bxi^{1/2}
\tbcW_t \bigr\rangle } \rright] %
,
\]
as $n \to\infty$.
This limiting random vector can be written in the form as given in Theorem~\ref{main_Ad}, since $\langle\bone, \bcX_t \rangle= \cY_t$,
$\langle\bone, \bcM_t \rangle
= \langle\bone, \bcX_t \rangle- \langle\bone, {\mathbf
{m}}_\bvare\rangle t
= \cY_t - \langle\bone, {\mathbf{m}}_\bvare\rangle t$
and
$\langle\tbu, \obV_\bxi^{1/2} \tbcW_t \rangle
= \langle\obV_\bxi\tbu, \tbu\rangle^{1/2} \tcW_t$
for all $t \in\RR_+$ with a (one-dimensional) standard Wiener process
$(\tcW_t)_{t\in\RR_+}$.
\begin{pf*}{Proof of Theorem~\ref{main_conv}}
In order to show convergence $\bcZ^{(n)} \distr\bcZ$, we apply Theorem~\ref{Conv2DiffThm} with the special choices $\bcU:= \bcZ$,
$\bU^{(n)}_k := \bZ^{(n)}_k$, $n, k \in\NN$,
$(\cF_k^{(n)})_{k\in\ZZ_+} := (\cF_k)_{k\in\ZZ_+}$ and the
function $\gamma$
which is defined in Theorem~\ref{main_conv}.
Note that the discussion after Theorem~\ref{main_conv} shows that the SDE
\eqref{ZSDE} admits a unique strong solution $(\bcZ_t^{\mathbf
{z}})_{t\in\RR_+}$ for
all initial values $\bcZ_0^{\mathbf{z}}= {\mathbf{z}}\in
(\RR^2)^3$.

Now we show that conditions (i) and (ii) of Theorem~\ref{Conv2DiffThm} hold.
The conditional variance
$\EE (\bZ^{(n)}_k (\bZ^{(n)}_k)^\top|\cF_{k-1} )$
has the form
\[
\lleft[\matrix{ n^{-2} & n^{-3}
U_{k-1} & n^{-5/2} V_{k-1}
\cr
n^{-3}
U_{k-1} & n^{-4} U_{k-1}^2 &
n^{-7/2} U_{k-1} V_{k-1}
\cr
n^{-5/2}
V_{k-1} & n^{-7/2} U_{k-1} V_{k-1} &
n^{-3} V_{k-1}^2 } \rright] %
\otimes\bV_{\bM_k}
\]
for $n \in\NN$, $k \in\{1, \ldots, n\}$, with
$\bV_{\bM_k} := \EE(\bM_k \bM_k^\top|\cF_{k-1})$, and
$\gamma(s,\bcZ_s^{(n)}) \gamma(s,\bcZ_s^{(n)})^\top$ has the form
\[
\lleft[\matrix{ \bigl\langle\bone, \bcM_s^{(n)}
+ s {\mathbf{m}}_\bvare\bigr\rangle & \bigl\langle\bone,
\bcM_s^{(n)} + s {\mathbf{m}}_\bvare\bigr
\rangle^2 & \bzero
\cr
\bigl\langle\bone, \bcM_s^{(n)}
+ s {\mathbf{m}}_\bvare\bigr\rangle^2 & \bigl\langle\bone,
\bcM_s^{(n)} + s {\mathbf{m}}_\bvare\bigr
\rangle^3 & \bzero
\cr
\bzero& \bzero & \displaystyle \frac{\langle\obV_\bxi\tbu, \tbu\rangle}{4 \alpha\beta} \bigl\langle
\bone, \bcM_s^{(n)} + s {\mathbf{m}}_\bvare\bigr
\rangle^2 } \rright] %
\otimes\obV_\xi
\]
for $s\in\RR_+$, where we used that
$\langle\bone, \bcM_s^{(n)} + s {\mathbf{m}}_\bvare\rangle^+
= \langle\bone, \bcM_s^{(n)} + s {\mathbf{m}}_\bvare\rangle$,
$s \in\RR_+$,
$n \in\NN$.
Indeed, by \eqref{Mk}, we get
%
%e5.4 #&#
\begin{eqnarray}
\label{M+} %
 \bigl\langle\bone,
\bcM_s^{(n)} + s {\mathbf{m}}_\bvare\bigr\rangle &=&
\frac{1}{n} \sum_{k=1}^{\lfloor ns\rfloor} \langle
\bone, \bX_k - {\mathbf{m}}_\bxi\bX_{k-1} - {
\mathbf{m}}_\bvare\rangle + \langle\bone, s {\mathbf{m}}_\bvare
\rangle\nonumber
\\
&=& \frac{1}{n} \langle\bone, \bX_{\lfloor ns\rfloor}\rangle +
\frac{ns - {\lfloor ns\rfloor}}{n} \langle\bone, {\mathbf {m}}_\bvare\rangle
\\
&=& \frac{1}{n} U_{\lfloor ns\rfloor}+ \frac{ns - {\lfloor ns\rfloor
}}{n} \langle\bone, {
\mathbf{m}}_\bvare\rangle \in\RR_+ \nonumber %
\end{eqnarray}
for $s \in\RR_+$, $n \in\NN$, since $\bone^\top{\mathbf
{m}}_\bxi= \bone^\top$
implies
$\langle\bone, {\mathbf{m}}_\bxi\bX_{k-1} \rangle
= \bone^\top{\mathbf{m}}_\bxi\bX_{k-1} = \bone^\top\bX_{k-1}
= \langle\bone, \bX_{k-1} \rangle$.

In order to check condition (i) of Theorem~\ref{Conv2DiffThm}, we need to
prove that for each $T>0$, as $n \to\infty$,
%
%e5.5 #&#
%e5.6 #&#
%e5.7 #&#
%e5.8 #&#
\begin{eqnarray}
\label{Zcond1}\sup_{t\in[0,T]} \Biggl\| \frac{1}{n^2} \sum
_{k=1}^{{\lfloor nt\rfloor}} \bV_{\bM_k} - \int
_0^t \bigl\langle\bone, \bcM_s^{(n)}
+ s {\mathbf{m}}_\bvare\bigr\rangle \obV_\bxi \,\mathrm{d}s \Biggr\|
&\stoch&0 ,
\\
\label{Zcond2}\sup_{t\in[0,T]} \Biggl\| \frac{1}{n^3} \sum
_{k=1}^{{\lfloor nt\rfloor}} U_{k-1} \bV_{\bM_k} -
\int_0^t \bigl\langle\bone,
\bcM_s^{(n)} + s {\mathbf{m}}_\bvare\bigr
\rangle^2 \obV_\bxi \,\mathrm{d}s \Biggr\| &\stoch&0 ,
\\
\label{Zcond3}\sup_{t\in[0,T]} \Biggl\| \frac{1}{n^4} \sum
_{k=1}^{{\lfloor nt\rfloor}} U_{k-1}^2
\bV_{\bM_k} - \int_0^t \bigl\langle
\bone, \bcM_s^{(n)} + s {\mathbf{m}}_\bvare\bigr
\rangle^3 \obV_\bxi \,\mathrm{d}s \Biggr\| &\stoch&0 ,
\\
\label{Zcond4}\sup_{t\in[0,T]} \Biggl\| \frac{1}{n^3} \sum
_{k=1}^{{\lfloor nt\rfloor}} V_{k-1}^2
\bV_{\bM_k} - \frac{\langle\obV_\bxi\tbu, \tbu\rangle}{4 \alpha\beta} \int_0^t
\bigl\langle\bone, \bcM_s^{(n)} + s {\mathbf{m}}_\bvare
\bigr\rangle^2 \obV_\bxi \,\mathrm{d}s \Biggr\| &\stoch&0 ,
\\
\label{Zcond5}\sup_{t\in[0,T]} \Biggl\| \frac{1}{n^{5/2}} \sum
_{k=1}^{{\lfloor nt\rfloor}} V_{k-1} \bV_{\bM_k} \Biggr\|
&\stoch&0 ,
\\
\label{Zcond6}\sup_{t\in[0,T]} \Biggl\| \frac{1}{n^{7/2}} \sum
_{k=1}^{{\lfloor nt\rfloor}} U_{k-1} V_{k-1}
\bV_{\bM_k} \Biggr\| &\stoch&0 .
\end{eqnarray}

First, we show \eqref{Zcond1}.
By \eqref{M+},
$\int_0^t \langle\bone, \bcM_s^{(n)} + s {\mathbf{m}}_\bvare
\rangle\,\mathrm{d}s$
has the form
\[
\frac{1}{n^2} \sum_{k=1}^{{\lfloor nt\rfloor}-1}
U_k + \frac{nt - {\lfloor nt\rfloor}}{n^2} U_{\lfloor nt\rfloor} + \frac{{\lfloor nt\rfloor}+ (nt - {\lfloor nt\rfloor})^2}{2 n^2}
\langle\bone, {\mathbf{m}}_\bvare\rangle.
\]
Using Lemma~\ref{Moments}, we obtain
%
%e5.11 #&#
\begin{equation}
\label{M2F} \bV_{\bM_k} = U_{k-1} \obV_\bxi+
\tfrac{1}{2} V_{k-1} (\bV_{\bxi_1} - \bV _{\bxi_2}) +
\bV_\bvare.
\end{equation}
Thus, in order to show \eqref{Zcond1}, it suffices to prove
%
%e5.12 #&#
%e5.13 #&#
\begin{eqnarray}
\label
{1supsumV_1supU}&n^{-2} \displaystyle \sum_{k=1}^{{\lfloor nT\rfloor}}
|V_k| \stoch0 ,\qquad  n^{-2} \displaystyle \sup_{t \in[0,T]}
U_{\lfloor nt\rfloor}\stoch0, &
\\
\label{1supnt}&n^{-2} \displaystyle \sup_{t \in[0,T]} \bigl[ {\lfloor nt\rfloor}+
\bigl(nt - {\lfloor nt\rfloor}\bigr)^2 \bigr] \to0 ,&
\end{eqnarray}
as $n \to\infty$.
Using \eqref{seged_UV_UNIFORM1} with $(\ell, i, j) = (2, 1, 1)$ and
\eqref{seged_UV_UNIFORM2} with $(\ell, i, j) = (2, 1, 0)$, we have
\eqref{1supsumV_1supU}.
Clearly, \eqref{1supnt} follows from $|nt - {\lfloor nt\rfloor}| \leq
1$, $n \in\NN$,
$t \in\RR_+$, thus we conclude \eqref{Zcond1}.
The convergences \eqref{Zcond2} and \eqref{Zcond3} can be checked in
a similar
way.

Next, we turn to prove \eqref{Zcond4}.
By \eqref{M2F} and \eqref{seged_UV_UNIFORM1}, we get
%
%e5.14 #&#
\begin{eqnarray}
\label{Zcond2a} n^{-3} \sup_{t \in[0,T]} \Biggl\llVert \sum
_{k=1}^{{\lfloor nt\rfloor}} U_{k-1}
\bV_{\bM_k} - \sum_{k=1}^{{\lfloor nt\rfloor}}
U_{k-1}^2 \obV_\bxi\Biggr\rrVert \stoch0 ,
\end{eqnarray}
as $n \to\infty$ for all $T>0$.
Using \eqref{Zcond2}, in order to prove \eqref{Zcond4}, it is
sufficient to
show that
%
%e5.15 #&#
\begin{eqnarray}
\label{Zcond4a} n^{-3} \sup_{t \in[0,T]} \Biggl\llVert \sum
_{k=1}^{{\lfloor nt\rfloor}} V_{k-1}^2
\bV_{\bM_k} - \frac{\langle\obV_\bxi\tbu, \tbu\rangle}{4 \alpha\beta} \sum_{k=1}^{{\lfloor nt\rfloor}}
U_{k-1}^2 \obV_\bxi\Biggr\rrVert \stoch0 ,
\end{eqnarray}
as $n \to\infty$ for all $T>0$.
By \eqref{M2F}, $\sum_{k=1}^{{\lfloor nt\rfloor}} V_{k-1}^2 \bV
_{\bM_k}$ has the form
\[
\sum_{k=1}^{{\lfloor nt\rfloor}} U_{k-1}
V_{k-1}^2 \obV_\bxi + \frac{1}{2} \sum
_{k=1}^{{\lfloor nt\rfloor}} V_{k-1}^3 (
\bV_{\bxi
_1} - \bV_{\bxi_2}) + \sum_{k=1}^{{\lfloor nt\rfloor}}
V_{k-1}^2 \bV_\bvare.
\]
Using \eqref{seged_UV_UNIFORM1} with $(\ell, i, j) = (6, 0, 3)$ and
$(\ell, i, j) = (4, 0, 2)$, we have
\begin{eqnarray*}
n^{-3} \sum_{k=1}^{{\lfloor nT\rfloor}} |
V_k |^3 \stoch0 ,\qquad  n^{-3} \sum
_{k=1}^{{\lfloor nT\rfloor}} V_k^2 \stoch0 ,
\qquad \mbox{as }n\to\infty,
\end{eqnarray*}
hence \eqref{Zcond4a} will follow from
%
%e5.16 #&#
\begin{eqnarray}
\label{Zcond4b} n^{-3} \sup_{t \in[0,T]} \Biggl\llVert \sum
_{k=1}^{{\lfloor nt\rfloor}} U_{k-1}
V_{k-1}^2 - \frac{\langle\obV_\bxi\tbu, \tbu\rangle}{4 \alpha\beta} \sum
_{k=1}^{{\lfloor nt\rfloor}} U_{k-1}^2 \Biggr
\rrVert \stoch0 ,
\end{eqnarray}
as $n \to\infty$ for all $T>0$.
By the method of the proof of Lemma~\ref{main_VV}, we obtain a decomposition
of $\sum_{k=1}^{{\lfloor nt\rfloor}} U_{k-1} V_{k-1}^2$ as a sum of a
martingale and some
negligible terms, namely,
\begin{eqnarray*}
\sum_{k=1}^{{\lfloor nt\rfloor}} U_{k-1}
V_{k-1}^2 &=& \frac{1}{4 \alpha\beta} \sum
_{k=2}^{{\lfloor nt\rfloor}} \bigl[U_{k-1}
V_{k-1}^2 - \EE\bigl(U_{k-1} V_{k-1}^2
|\cF_{k-2}\bigr) \bigr]
\\
&&{}+ \frac{\langle\obV_\bxi\tbu, \tbu\rangle}{4 \alpha\beta} \sum_{k=2}^{{\lfloor nt\rfloor}}
U_{k-2}^2 - \frac{(\alpha- \beta)^2}{4 \alpha\beta} U_{{\lfloor nt\rfloor}-
1}
V_{{\lfloor nt\rfloor}- 1}^2 + \OO(n)
\\
&&{}+ \mbox{lin. comb. of }\sum_{k=2}^{{\lfloor nt\rfloor}}
U_{k-2} V_{k-2}, \sum_{k=2}^{{\lfloor nt\rfloor}}
V_{k-2}^2, \sum_{k=2}^{{\lfloor
nt\rfloor}}
U_{k-2}\mbox{ and }\sum_{k=2}^{{\lfloor nt\rfloor}}
V_{k-2}. %
\end{eqnarray*}
Using \eqref{seged_UV_UNIFORM4} with $(\ell, i, j) = (8, 1, 2)$ we have
\[
n^{-3}\sup_{t \in[0,T]} \Biggl\vert\sum
_{k=2}^{\lfloor nt\rfloor} \bigl[U_{k-1}
V_{k-1}^2 - \EE\bigl(U_{k-1} V_{k-1}^2
\Biggr|\cF_{k-2}\bigr) \bigr] \vert \stoch0 , \qquad \mbox{as }n\to\infty.
\]
Thus, in order to show \eqref{Zcond4b}, it suffices to prove
%
%e5.17 #&#
%e5.18 #&#
%e5.19 #&#
\begin{eqnarray}
\label{4supsumUV_4supsumVV}n^{-3} \sum_{k=1}^{{\lfloor nT\rfloor}} |
U_k V_k | &\stoch&0 ,\qquad  n^{-3} \sum
_{k=1}^{{\lfloor nT\rfloor}} V_k^2 \stoch0 ,
\\
\label{4supsumU_4supsumV}n^{-3} \sum_{k=1}^{{\lfloor nT\rfloor}}
U_k &\stoch&0 ,\qquad  n^{-3} \sum_{k=1}^{{\lfloor nT\rfloor}}
|V_k| \stoch0 ,
\\
\label{4supUVV_4supUU}n^{-3} \sup_{t \in[0,T]} U_{\lfloor nt\rfloor}V_{\lfloor nt\rfloor
}^2
&\stoch&0 ,\qquad  n^{-3/2} \sup_{t \in[0,T]} U_{{\lfloor nt\rfloor}} \stoch0
,
\end{eqnarray}
as $n \to\infty$.
Using \eqref{seged_UV_UNIFORM1} with $(\ell, i, j) = (2, 1, 1)$,
$(\ell, i, j) = (4, 0, 2)$, $(\ell, i, j) = (2, 1, 0)$ and
$(\ell, i, j) = (2, 0, 1)$, we have \eqref{4supsumUV_4supsumVV} and
\eqref{4supsumU_4supsumV}.
By \eqref{seged_UV_UNIFORM2} with $(\ell, i, j) = (4, 1, 2)$ and by
\eqref{seged_UV_UNIFORM2}, we have \eqref{4supUVV_4supUU}.
Thus, we conclude \eqref{Zcond4}.
Convergences \eqref{Zcond5} and \eqref{Zcond6} can be proved similarly.

Finally, we check condition (ii) of Theorem~\ref{Conv2DiffThm}, that
is, the
conditional Lindeberg condition
%
%e5.20 #&#
\begin{equation}
\label{Zcond3_new} \sum_{k=1}^{\lfloor nT \rfloor} \EE \bigl( \bigl\|
\bZ^{(n)}_k\bigr\|^2 \bbone_{\{\|\bZ^{(n)}_k\| > \theta\}} |
\cF_{k-1} \bigr) \stoch0 ,\qquad  \mbox{as }n\to\infty
\end{equation}
for all $\theta>0$ and $T>0$.
We have
$\EE ( \|\bZ^{(n)}_k\|^2 \bbone_{\{\|\bZ^{(n)}_k\| > \theta\}}
|\cF_{k-1}  )
\leq\theta^{-2} \EE ( \|\bZ^{(n)}_k\|^4  |\cF_{k-1}  )$
and
\[
\bigl\|\bZ^{(n)}_k\bigr\|^4 \leq3 \bigl(
n^{-4} + n^{-8} U_{k-1}^4 +
n^{-6} V_{k-1}^4 \bigr) \|\bM_{k-1}
\|^4 .
\]
Hence, for all $\theta>0$ and $T>0$, we have
\[
\sum_{k=1}^{{\lfloor nT\rfloor}} \EE \bigl( \bigl\|
\bZ^{(n)}_k\bigr\|^2 \bbone_{\{\|\bZ^{(n)}_k\| > \theta\}} \bigr)
\to0 , \qquad \mbox{as }n\to\infty,
\]
since $\EE( \|\bM_k\|^4 ) = \OO(k^2)$,
$\EE( \|\bM_k\|^4 U_{k-1}^4 ) \leq\sqrt{\EE(\|\bM_k\|^8) \EE(U_{k-1}^8)}
= \OO(k^6)$
and\linebreak[4]
$\EE( \|\bM_k\|^4 V_{k-1}^4 ) \leq\sqrt{\EE(\|\bM_k\|^8) \EE(V_{k-1}^8)}
= \OO(k^4)$
by Corollary~\ref{EEX_EEU_EEV}.
Here we call the attention that our eighth order moment conditions
$\EE(\|\bxi_{1,1,1}\|^8) < \infty$, $\EE(\|\bxi_{1,1,2}\|^8) <
\infty$ and
$\EE(\|\bvare_1\|^8) < \infty$ are used for applying Corollary~\ref{EEX_EEU_EEV}.
This yields \eqref{Zcond3_new}.
\end{pf*}

%s6 #&#
\section{Proof of Theorem \texorpdfstring{\protect\ref{main1_Ad}}{4.2}}
\label{section_proof_main1}

This is similar to the proof of Theorem~\ref{main_Ad}.
Consider the sequence of stochastic processes
\[
\bcZ^{(n)}_t := %
\lleft[\matrix{
\bcM_t^{(n)}
\cr
\cN_t^{(n)}
\cr
\bcP_t^{(n)} } \rright] %
:= \sum
_{k=1}^{\lfloor nt\rfloor} \bZ^{(n)}_k \qquad \mbox{with } \bZ^{(n)}_k := %
\lleft[\matrix{
n^{-1} \bM_k
\cr
n^{-3/2} \langle\bone,
\bM_k \rangle U_{k-1}
\cr
n^{-3/2} \bM_k
V_{k-1} } \rright] %
\]
for $t \in\RR_+$ and $k, n \in\NN$.
Theorem~\ref{main1_Ad} follows from Lemma~\ref{main_VV} and the following
theorem (this will be explained after Theorem~\ref{main1_conv}).
%
%th6.1 #&#
\begin{Thm}\label{main1_conv}
If $\langle\obV_\bxi\bone, \bone\rangle= 0$ then
%
%e6.1 #&#
\begin{equation}
\label{conv_Z1} \bcZ^{(n)} \distr\bcZ, \qquad \mbox{as }n\to\infty,
\end{equation}
where the process $(\bcZ_t)_{t \in\RR_+}$ with values in
$\RR^2 \times\RR\times\RR^2$ is the unique strong solution of the SDE
%
%e6.2 #&#
\begin{equation}
\label{ZSDE1} \mathrm{d}\bcZ_t = \gamma(t, \bcZ_t)
\lleft[\matrix{ \mathrm{d}\bcW_t
\cr
\mathrm{d}
\ttcW_t
\cr
\mathrm {d}\tbcW_t } \rright] %
, \qquad t
\in\RR_+ ,
\end{equation}
with initial value $\bcZ_0 = \bzero$, where $(\bcW_t)_{t \in\RR_+}$,
$(\ttcW_t)_{t \in\RR_+}$ and $(\tbcW_t)_{t \in\RR_+}$ are independent
standard Wiener processes of dimension $2$, $1$ and $2$, respectively, and
$\gamma(t, {\mathbf{x}})$ is a block diagonal matrix with the matrices
$\langle\bone, ({\mathbf{x}}_1 + t {\mathbf{m}}_\bvare)^+
\rangle^{1/2} \obV_\bxi^{1/2}$,
$\langle\bV_\bvare\bone, \bone\rangle^{1/2}
\langle\bone, {\mathbf{m}}_\bvare\rangle t$
and
$ (\frac{\langle\obV_\bxi\tbu, \tbu\rangle} {
4 \alpha\beta} )^{1/2}
\langle\bone, {\mathbf{x}}_1 + t {\mathbf{m}}_\bvare
\rangle\obV_\bxi^{1/2}$
in its diagonal for each
$t \in\RR_+$ and
${\mathbf{x}}= ({\mathbf{x}}_1 , x_2 , {\mathbf{x}}_3)
\in\RR^2 \times\RR\times\RR^2$.
\end{Thm}

As in the case of Theorem~\ref{main_Ad}, the SDE \eqref{ZSDE1} has a unique
strong solution with initial value $\bcZ_0 = \bzero$, for which we have
\[
\bcZ_t = %
\lleft[\matrix{ \bcM_t
\cr
\cN_t
\cr
\bcP_t } \rright] %
= %
\lleft[\matrix{ \displaystyle \int_0^t
\cY_t^{1/2} \obV_\bxi^{1/2} \,\mathrm{d}
\bcW_s
\cr
\langle\bV_\bvare\bone, \bone
\rangle^{1/2} \langle\bone, {\mathbf{m}}_\bvare\rangle \displaystyle \int
_0^t s \,\mathrm{d}\ttcW_s
\cr
\biggl(\displaystyle \frac{\langle\obV_\bxi\tbu, \tbu\rangle} {
4 \alpha\beta} \biggr)^{1/2} \int_0^t
\cY_t \obV_\bxi^{1/2} \,\mathrm{d}
\tbcW_s } \rright] %
, \qquad t\in\RR_+ ,
\]
where now $\langle\obV_\bxi\bone, \bone\rangle= 0$ yields
$\cY_t = \langle\bone, {\mathbf{m}}_\bvare\rangle t$, $t \in
\RR_+$.
One can again easily derive
%
%e6.3 #&#
\begin{eqnarray}
\label{convXZ1} %
\lleft[\matrix{ \bcX^{(n)}
\cr
\bcZ^{(n)} } \rright] %
\distr %
\lleft[\matrix{
\bcX
\cr
\bcZ } \rright] %
,\qquad  \mbox{as }n \to\infty,
\end{eqnarray}
where
\[
\bcX^{(n)}_t := n^{-1} \bX_{\lfloor nt\rfloor},\qquad
\bcX_t: = \frac{1}{2} \langle\bone, \bcM_t + t {
\mathbf {m}}_\bvare\rangle\bone = \frac{t}{2} \langle\bone, {
\mathbf{m}}_\bvare\rangle\bone,
\]
for $t \in\RR_+$ and $n\in\NN$, since
$\bcX_t = \frac{1}{2} \cY_t \bone
= \frac{t}{2} \langle\bone, {\mathbf{m}}_\bvare\rangle\bone$,
$t \in\RR_+$.
Next, similarly to the proof of \eqref{seged2}, by Lemma~\ref{Marci},
convergence \eqref{convXZ1} and Lemma~\ref{main_VV} with
$U_{k-1} = \langle\bone, \bX_{k-1} \rangle$ imply
\[
\sum_{k=1}^n %
\lleft[
\matrix{ n^{-3} U_{k-1}^2
\cr
n^{-2}
V_{k-1}^2
\cr
n^{-3/2} \langle\bone,
\bM_k \rangle U_{k-1}
\cr
n^{-3/2} \langle\tbu,
\bM_k \rangle V_{k-1} } \rright] %
\distr
\lleft[\matrix{ \displaystyle \int_0^1
\langle\bone, \bcX_t \rangle^2 \,\mathrm{d}t
\cr
\biggl(
\displaystyle \frac{\langle\obV_\bxi\tbu, \tbu\rangle} {
4 \alpha\beta} \biggr) \int_0^1 \langle
\bone, \bcX_t \rangle\,\mathrm{d}t
\cr
\langle\bV_\bvare\bone,
\bone\rangle^{1/2} \langle\bone, {\mathbf{m}}_\bvare\rangle \displaystyle \int
_0^1 t \,\mathrm{d}\ttcW_t
\cr
\biggl(\displaystyle \frac{\langle\obV_\bxi\tbu, \tbu\rangle} {
4 \alpha\beta} \biggr)^{1/2} \int_0^1
\cY_t \,\mathrm{d}\bigl\langle\tbu, \obV_\bxi^{1/2}
\tbcW_t \bigr\rangle } \rright] %
,
\]
as $n \to\infty$.
This limiting random vector can be written in the form as given in Theorem~\ref{main1_Ad} since
$\langle\bone, \bcX_t \rangle= \cY_t
= \langle\bone, {\mathbf{m}}_\bvare\rangle t$,
and
$\langle\tbu, \obV_\bxi^{1/2} \tbcW_t \rangle
= \langle\obV_\bxi\tbu, \tbu\rangle^{1/2} \tcW_t$
for all $t \in\RR_+$ with a\vspace*{1.5pt} (one-dimensional) standard Wiener process
$(\tcW_t)_{t\in\RR_+}$.
\begin{pf*}{Proof of Theorem~\ref{main1_conv}}
Similar to the proof of Theorem~\ref{main_conv}.
The conditional variance
$\EE (\bZ^{(n)}_k (\bZ^{(n)}_k)^\top|\cF_{k-1} )$
has the form
\[
\lleft[\matrix{ n^{-2} \bV_{\bM_k} &
n^{-5/2} U_{k-1} \bV_{\bM_k} \bone & n^{-5/2}
V_{k-1} \bV_{\bM_k}
\vspace*{2pt}\cr
n^{-5/2} U_{k-1}
\bone^\top\bV_{\bM_k} & n^{-3} U_{k-1}^2
\bone^\top\bV_{\bM_k} \bone & n^{-3} U_{k-1}
V_{k-1} \bone^\top\bV_{\bM_k}
\cr
n^{-5/2}
V_{k-1} \bV_{\bM_k} & n^{-3} U_{k-1}
V_{k-1} \bV_{\bM_k} \bone & n^{-3}
V_{k-1}^2 \bV_{\bM_k} } \rright] %
\]
for $n \in\NN$, $k \in\{1, \ldots, n\}$, with
$\bV_{\bM_k} := \EE( \bM_k \bM_k^\top|\cF_{k-1} )$, and
$\gamma(s,\bcZ_s^{(n)}) \gamma(s,\bcZ_s^{(n)})^\top$ has the form
\[
\lleft[\matrix{ \bigl\langle\bone, \bcM_s^{(n)}
+ s {\mathbf{m}}_\bvare\bigr\rangle\obV _\xi& \bzero &
\bzero
\cr
\bzero & \langle\bV_\bvare\bone, \bone\rangle \langle\bone, {
\mathbf{m}}_\bvare\rangle^2 s^2 & \bzero
\cr
\bzero& \bzero & \displaystyle \frac{\langle\obV_\bxi\tbu, \tbu\rangle}{4 \alpha\beta} \bigl\langle\bone, \bcM_s^{(n)}
+ s {\mathbf{m}}_\bvare\bigr\rangle^2 \obV_\xi }
\rright] %
\]
for $s\in\RR_+$.

In order to check condition (i) of Theorem~\ref{Conv2DiffThm}, we need to
prove only that for each $T>0$,
%
%e6.4 #&#
%e6.5 #&#
%e6.6 #&#
\begin{eqnarray}
\label{1Zcond2}\sup_{t\in[0,T]} \Biggl\| \frac{1}{n^{5/2}} \sum
_{k=1}^{{\lfloor nt\rfloor}} U_{k-1} \bone^\top
\bV_{\bM_k} \Biggr\| &\stoch&0 ,
\\
\label{1Zcond3}\sup_{t\in[0,T]} \Biggl| \frac{1}{n^3} \sum
_{k=1}^{{\lfloor nt\rfloor}} U_{k-1}^2
\bone^\top\bV_{\bM_k} \bone - \int_0^t
\langle\bV_\bvare\bone, \bone\rangle \langle\bone, {
\mathbf{m}}_\bvare\rangle^2 s^2 \,\mathrm{d}s \Biggr|
&\stoch&0 ,
\\
\label{1Zcond6}\sup_{t\in[0,T]} \Biggl\| \frac{1}{n^3} \sum
_{k=1}^{{\lfloor nt\rfloor}} U_{k-1} V_{k-1}
\bone^\top\bV_{\bM_k} \Biggr\| &\stoch&0 ,
\end{eqnarray}
as $n \to\infty$, since the rest, namely, \eqref{Zcond1}, \eqref{Zcond4}
and \eqref{Zcond5} have already been proved.

Clearly, $\langle\obV_\bxi\bone, \bone\rangle= 0$ implies
$\langle\bV_{\bxi_1} \bone, \bone\rangle= 0$ and
$\langle\bV_{\bxi_2} \bone, \bone\rangle= 0$.
For each $i \in\{1, 2\}$, we have
$\langle\bV_{\bxi_i} \bone, \bone\rangle= \bone^\top\bV_{\bxi
_i} \bone
= (\bV_{\bxi_i}^{1/2} \bone)^\top(\bV_{\bxi_i}^{1/2} \bone)
= \|\bV_{\bxi_i}^{1/2} \bone\|^2$,
hence we obtain $\bV_{\bxi_i}^{1/2} \bone= \bzero$, thus
$\bV_{\bxi_i} \bone= \bV_{\bxi_i}^{1/2} (\bV_{\bxi_i}^{1/2}
\bone) = \bzero$,
and hence $\bone^\top\bV_{\bxi_i} = \bzero$, implying also
$\bone^\top\obV_\bxi= \bzero$.

First we show \eqref{1Zcond2}.
By \eqref{M2F}, $\bone^\top\obV_\bxi= \bzero$ and
$\bone^\top\bV_{\bxi_i} = \bzero$ for $i \in\{1, 2\}$, we obtain
%
%e6.7 #&#
\begin{equation}
\label{1Zcond2m} \sum_{k=1}^{{\lfloor nt\rfloor}}
U_{k-1} \bone^\top\bV_{\bM_k} = \sum
_{k=1}^{{\lfloor nt\rfloor}} U_{k-1} \bone^\top
\bV_\bvare,
\end{equation}
hence using \eqref{seged_UV_UNIFORM1} with $(\ell, i, j) = (2, 1,
0)$, we
conclude \eqref{1Zcond2}.

Now we turn to check \eqref{1Zcond3}.
By \eqref{M2F},
\[
\sum_{k=1}^{{\lfloor nt\rfloor}} U_{k-1}^2
\bone^\top\bV_{\bM_k} \bone = \sum
_{k=1}^{{\lfloor nt\rfloor}} U_{k-1}^2
\bone^\top\bV_\bvare \bone = \sum_{k=1}^{{\lfloor nt\rfloor}}
U_{k-1}^2 \langle\bV_\bvare \bone, \bone\rangle,
\]
hence, in order to show \eqref{1Zcond3}, it suffices to prove
%
%e6.8 #&#
\begin{equation}
\label{1Zcond3m} \sup_{t\in[0,T]} \Biggl| \frac{1}{n^3} \sum
_{k=1}^{{\lfloor nt\rfloor}} U_{k-1}^2 -
\frac{t^3}{3} \langle\bone, {\mathbf{m}}_\bvare\rangle^2
\Biggr| \stoch0 ,\qquad  \mbox{as }n \to\infty.
\end{equation}
We have
\begin{eqnarray*}
\Biggl| \frac{1}{n^3} \sum_{k=1}^{{\lfloor nt\rfloor}}
U_{k-1}^2 - \frac{t^3}{3} \langle\bone, {
\mathbf{m}}_\bvare\rangle^2 \Biggr| &\leq&\frac{1}{n^3} \sum
_{k=1}^{{\lfloor nt\rfloor}} \bigl\llvert U_{k-1}^2
- (k-1)^2 \langle\bone, {\mathbf{m}}_\bvare
\rangle^2 \bigr\rrvert
\\
&&{} + \Biggl| \frac{1}{n^3} \sum_{k=1}^{{\lfloor nt\rfloor}}
(k-1)^2 - \frac{t^3}{3} \Biggr| \langle\bone, {\mathbf{m}}_\bvare
\rangle^2 ,
\end{eqnarray*}
where
\[
\sup_{t\in[0,T]} \Biggl| \frac{1}{n^3} \sum
_{k=1}^{{\lfloor nt\rfloor}} (k-1)^2 - \frac{t^3}{3} \Biggr|
\to0 , \qquad \mbox{as }n \to\infty,
\]
hence, in order to show \eqref{1Zcond3}, it suffices to prove
%
%e6.9 #&#
\begin{equation}
\label{1Zcond3mm} \frac{1}{n^3} \sum_{k=1}^{{\lfloor nT\rfloor}}
\bigl\llvert U_k^2 - k^2 \langle\bone, {
\mathbf {m}}_\bvare\rangle^2 \bigr\rrvert \stoch0 ,
\qquad \mbox{as }n \to\infty.
\end{equation}
For all $k \in\NN$, by Remark~\ref{REMARK0},
$\langle\obV_\bxi\bone, \bone\rangle= 0$ implies
\begin{eqnarray*}
U_k &=& \sum_{j=1}^{X_{k-1,1}} (
\xi_{k,j,1,1} + \xi_{k,j,1,2}) + \sum_{j=1}^{X_{k-1,2}}
(\xi_{k,j,2,1} + \xi_{k,j,2,2}) + (\vare_{k,1} +
\vare_{k,2})
\\
&\ase& X_{k-1,1} + X_{k-1,2} + \vare_{k,1} +
\vare_{k,2} = U_{k-1} + \langle\bone, \bvare_k
\rangle,
\end{eqnarray*}
hence $U_k = \sum_{i=1}^k \langle\bone, \bvare_i \rangle$.
By Kolmogorov's maximal inequality,\vspace*{-1pt}
\begin{eqnarray*}
\PP \Bigl( n^{-1} \max_{k\in\{1,\ldots,{\lfloor nT\rfloor}\}} \bigl| U_k - k
\langle\bone, {\mathbf{m}}_\bvare\rangle\bigr| \geq\vare \Bigr) &\leq&
n^{-2} \vare^{-2} \var(U_{\lfloor nT\rfloor})
\\
&=& \frac{{\lfloor nT\rfloor}}{n^2 \vare^2} \var\bigl(\langle\bone, \bvare_1
\rangle^2\bigr) \to0
\end{eqnarray*}
as $n \to\infty$ for all $\vare> 0$, thus\vspace*{-1pt}
\[
n^{-1} \max_{k\in\{1,\ldots,{\lfloor nT\rfloor}\}} \bigl| U_k - k \langle
\bone, {\mathbf{m}}_\bvare\rangle\bigr| \stoch0 , \qquad \mbox{as }n \to\infty.
\]
We have\vspace*{-1pt}
\[
\bigl| U_k^2 - k^2 \langle\bone, {
\mathbf{m}}_\bvare\rangle^2 \bigr| \leq\bigl| U_k - k
\langle\bone, {\mathbf{m}}_\bvare\rangle\bigr|^2 + 2 k \langle
\bone, {\mathbf{m}}_\bvare\rangle \bigl| U_k - k \langle\bone, {
\mathbf{m}}_\bvare\rangle\bigr| ,
\]
hence
\begin{eqnarray*}
n^{-2} \max_{k\in\{1,\ldots,{\lfloor nT\rfloor}\}} \bigl| U_k^2 -
k^2 \langle\bone, {\mathbf{m}}_\bvare\rangle^2 \bigr|
&\leq& \Bigl( n^{-1} \max_{k\in\{1,\ldots,{\lfloor nT\rfloor}\}} \bigl| U_k - k
\langle\bone, {\mathbf{m}}_\bvare\rangle\bigr| \Bigr)^2
\\
&&{}+ \frac{2 {\lfloor nT\rfloor}}{n^2} \langle\bone, {\mathbf {m}}_\bvare\rangle \max
_{k\in\{1,\ldots,{\lfloor nT\rfloor}\}} \bigl| U_k - k \langle\bone, {
\mathbf{m}}_\bvare\rangle\bigr| \stoch0 ,
\end{eqnarray*}
as $n \to\infty$.
Consequently,
\begin{eqnarray*}
&&\frac{1}{n^3} \sum_{k=1}^{{\lfloor nT\rfloor}} \bigl
\llvert U_{k-1}^2 - (k-1)^2 \langle\bone, {
\mathbf{m}}_\bvare \rangle^2 \bigr\rrvert
\\
&&\quad \leq\frac{{\lfloor nT\rfloor}}{n^3} \max_{k\in\{1,\ldots,{\lfloor nT\rfloor}\}} \bigl\llvert
U_{k-1}^2 - (k-1)^2 \langle\bone, {
\mathbf{m}}_\bvare \rangle^2 \bigr\rrvert \stoch0 ,
\end{eqnarray*}
as $n \to\infty$, thus we conclude \eqref{1Zcond3mm}, and hence
\eqref{1Zcond3}.

Finally, we check \eqref{1Zcond6}.
By \eqref{M2F},
\[
\sum_{k=1}^{{\lfloor nt\rfloor}} U_{k-1}
V_{k-1} \bone^\top\bV _{\bM_k} = \sum
_{k=1}^{{\lfloor nt\rfloor}} U_{k-1} V_{k-1}
\bone^\top\bV _\bvare,
\]
hence using \eqref{seged_UV_UNIFORM1} with $(\ell, i, j) = (2, 1,
1)$, we
conclude \eqref{1Zcond6}.
Condition (ii) of Theorem~\ref{Conv2DiffThm} can be checked as in case of
Theorem~\ref{main_conv}.
\end{pf*}

%s7 #&#
\section{Proof of Theorem \texorpdfstring{\protect\ref{maint_Ad}}{4.3}}
\label{section_proof_maint}

This proof is also similar to the proof of Theorem~\ref{main_Ad}.
Consider the sequence of stochastic processes
\[
\bcZ^{(n)}_t := %
\lleft[\matrix{
\bcM_t^{(n)}
\cr
\bcN_t^{(n)}
\cr
\cP_t^{(n)} } \rright] %
:= \sum
_{k=1}^{\lfloor nt\rfloor} \bZ^{(n)}_k \qquad \mbox{with } \bZ^{(n)}_k := %
\lleft[\matrix{
n^{-1} \bM_k
\cr
n^{-2} \bM_k
U_{k-1}
\cr
n^{-1/2} \langle\tbu, \bM_k \rangle
V_{k-1} } \rright] %
\]
for $t \in\RR_+$ and $k, n \in\NN$.
Theorem~\ref{maint_Ad} follows from Lemma~\ref{main_VVt} and the following
theorem (this will be explained after Theorem~\ref{maint_conv}).
%
%th7.1 #&#
\begin{Thm}\label{maint_conv}
If $\langle\obV_\bxi\tbu, \tbu\rangle= 0$ then
%
%e7.1 #&#
\begin{equation}
\label{conv_Zt} \bcZ^{(n)} \distr\bcZ, \qquad \mbox{as }n\to\infty,
\end{equation}
where the process $(\bcZ_t)_{t \in\RR_+}$ with values in
$\RR^2 \times\RR^2 \times\RR$ is the unique strong solution of the SDE
%
%e7.2 #&#
\begin{equation}
\label{ZSDEt} \mathrm{d}\bcZ_t = \gamma(t, \bcZ_t)
\lleft[\matrix{ \mathrm{d}\bcW_t
\cr
\mathrm{d}
\tcW_t } \rright] %
, \qquad t \in\RR_+ ,
\end{equation}
with initial value $\bcZ_0 = \bzero$, where $(\bcW_t)_{t \in\RR
_+}$ and
$(\tcW_t)_{t \in\RR_+}$ are independent standard Wiener processes of
dimension $2$ and $1$, respectively, and
$\gamma\dvtx  \RR_+ \times(\RR^2 \times\RR^2 \times\RR) \to\RR
^{5\times3}$
is defined by
\[
\gamma(t, {\mathbf{x}}) := %
\lleft[\matrix{ \bigl\langle\bone,
({\mathbf{x}}_1 + t {\mathbf{m}}_\bvare)^+ \bigr
\rangle^{1/2} \obV_\bxi^{1/2} & \bzero
\cr
\bigl
\langle\bone, ({\mathbf{x}}_1 + t {\mathbf{m}}_\bvare)^+
\bigr\rangle^{3/2} \obV_\bxi^{1/2} & \bzero
\cr
\bzero& \bigl[\langle\bV_\bvare\tbu, \tbu\rangle \EE\bigl(\langle\tbu,
\bvare_1 \rangle^2\bigr) \bigr]^{1/2} } \rright] %
\]
for $t \in\RR_+$ and
${\mathbf{x}}= ({\mathbf{x}}_1 , {\mathbf{x}}_2 , x_3)
\in\RR^2 \times\RR^2 \times\RR$.
\end{Thm}

As in the case of Theorem~\ref{main_Ad}, the SDE \eqref{ZSDEt} has a unique
strong solution with initial value $\bcZ_0 = \bzero$, for which we have
\[
\bcZ_t = %
\lleft[\matrix{ \bcM_t
\cr
\bcN_t
\cr
\cP_t } \rright] %
= %
\lleft[\matrix{ \displaystyle \int_0^t
\cY_s^{1/2} \obV_\bxi^{1/2} \,\mathrm{d}
\bcW_s
\cr
\displaystyle \int_0^t
\cY_s \,\mathrm{d}\bcM_s
\cr
\bigl[\langle
\bV_\bvare\tbu, \tbu\rangle \EE\bigl(\langle\tbu, \bvare_1
\rangle^2\bigr) \bigr]^{1/2} \tcW_t } \rright] %
,\qquad  t\in\RR_+ .
\]
One can again easily derive
%
%e7.3 #&#
\begin{eqnarray}
\label{convXZt} %
\lleft[\matrix{ \bcX^{(n)}
\cr
\bcZ^{(n)} } \rright] %
\distr %
\lleft[\matrix{
\bcX
\cr
\bcZ } \rright] %
, \qquad \mbox{as }n \to\infty,
\end{eqnarray}
where
\[
\bcX^{(n)}_t := n^{-1} \bX_{\lfloor nt\rfloor},\qquad
\bcX_t: = \tfrac{1}{2} \langle\bone, \bcM_t + t {
\mathbf {m}}_\bvare\rangle\bone,\qquad  t \in\RR_+ ,\ n \in\NN.
\]
Next, similarly to the proof of \eqref{seged2}, by Lemma~\ref{Marci},
convergence \eqref{convXZt} and Lemma~\ref{main_VVt} imply
\[
\sum_{k=1}^n %
\lleft[
\matrix{ n^{-3} U_{k-1}^2
\cr
n^{-1}
V_{k-1}^2
\cr
n^{-2} \langle\bone,
\bM_k \rangle U_{k-1}
\cr
n^{-1/2} \langle\tbu,
\bM_k \rangle V_{k-1} } \rright] %
\distr
\lleft[\matrix{ \displaystyle \int_0^1
\langle\bone, \bcX_t \rangle^2 \,\mathrm{d}t
\cr
\EE\bigl(
\langle\tbu, \bvare_1 \rangle^2\bigr)
\cr
\displaystyle \int
_0^1 \cY_t \,\mathrm{d}\langle\bone,
\bcM_t \rangle
\cr
\bigl[\langle\bV_\bvare\tbu, \tbu\rangle
\EE\bigl(\langle\tbu, \bvare_1 \rangle^2\bigr)
\bigr]^{1/2} \tcW_1 } \rright] %
,
\]
as $n \to\infty$.
Note that this convergence holds even in case
$\EE[\langle\tbu, \bvare_1 \rangle^2] = 0$.
The limiting random vector can be written in the form as given in Theorem~\ref{maint_Ad}, since $\langle\bone, \bcX_t \rangle= \cY_t$ and
$\langle\bone, \bcM_t \rangle
= \cY_t - \langle\bone, {\mathbf{m}}_\bvare t \rangle$
for all $t \in\RR_+$.
\begin{pf*}{Proof of Theorem~\ref{maint_conv}}
Similar to the proof of Theorem~\ref{main_conv}.
The conditional variance
$\EE (\bZ^{(n)}_k (\bZ^{(n)}_k)^\top|\cF_{k-1} )$
has the form
\[
\lleft[\matrix{ n^{-2} \bV_{\bM_k} &
n^{-3} U_{k-1} \bV_{\bM_k} & n^{-3/2}
V_{k-1} \bV_{\bM_k} \tbu
\cr
n^{-3} U_{k-1}
\bV_{\bM_k} & n^{-4} U_{k-1}^2
\bV_{\bM_k} & n^{-5/2} U_{k-1} V_{k-1}
\bV_{\bM_k} \tbu
\cr
n^{-3/2} V_{k-1} \tbu^\top
\bV_{\bM_k} & n^{-5/2} U_{k-1} V_{k-1}
\tbu^\top\bV_{\bM_k} & n^{-1} V_{k-1}^2
\tbu^\top\bV_{\bM_k} \tbu } \rright] %
\]
for $n \in\NN$, $k \in\{1, \ldots, n\}$, with
$\bV_{\bM_k} := \EE( \bM_k \bM_k^\top|\cF_{k-1} )$, and
$\gamma(s,\bcZ_s^{(n)}) \gamma(s,\bcZ_s^{(n)})^\top$ has the form
\[
\lleft[\matrix{ \bigl\langle\bone, \bcM_s^{(n)}
+ s {\mathbf{m}}_\bvare\bigr\rangle\obV _\xi & \bigl\langle
\bone, \bcM_s^{(n)} + s {\mathbf{m}}_\bvare\bigr
\rangle^2 \obV_\xi & \bzero
\cr
\bigl\langle\bone,
\bcM_s^{(n)} + s {\mathbf{m}}_\bvare\bigr
\rangle^2 \obV_\xi & \bigl\langle\bone,
\bcM_s^{(n)} + s {\mathbf{m}}_\bvare\bigr
\rangle^3 \obV_\xi & \bzero
\cr
\bzero& \bzero & \langle
\bV_\bvare\tbu, \tbu\rangle \EE\bigl(\langle\tbu, \bvare_1
\rangle^2\bigr) } \rright] %
\]
for $s\in\RR_+$.

In order to check condition (i) of Theorem~\ref{Conv2DiffThm}, we need to
prove only that for each $T>0$,
%
%e7.4 #&#
%e7.5 #&#
%e7.6 #&#
\begin{eqnarray}
\label{tZcond4}\sup_{t\in[0,T]} \Biggl| \frac{1}{n} \sum
_{k=1}^{{\lfloor nt\rfloor}} V_{k-1}^2
\tbu^\top\bV_{\bM_k} \tbu - t \langle\bV_\bvare\tbu,
\tbu\rangle \EE\bigl(\langle\tbu, \bvare_1 \rangle^2\bigr)
\Biggr| &\stoch&0 ,
\\
\label{tZcond5}\sup_{t\in[0,T]} \Biggl\| \frac{1}{n^{3/2}} \sum
_{k=1}^{{\lfloor nt\rfloor}} V_{k-1} \tbu^\top
\bV_{\bM_k} \Biggr\| &\stoch&0 ,
\\
\label{tZcond6}\sup_{t\in[0,T]} \Biggl\| \frac{1}{n^{5/2}} \sum
_{k=1}^{{\lfloor nt\rfloor}} U_{k-1} V_{k-1}
\tbu^\top\bV_{\bM_k} \Biggr\| &\stoch&0 ,
\end{eqnarray}
as $n \to\infty$, since the rest, namely, \eqref{Zcond1}, \eqref{Zcond2}
and \eqref{Zcond3}, have already been proved.

Clearly, $\langle\obV_\bxi\tbu, \tbu\rangle= 0$ implies
$\langle\bV_{\bxi_1} \tbu, \tbu\rangle= 0$ and
$\langle\bV_{\bxi_2} \tbu, \tbu\rangle= 0$.
For each $i \in\{1, 2\}$, we have
$\langle\bV_{\bxi_i} \tbu, \tbu\rangle= \tbu^\top\bV_{\bxi_i}
\tbu
= (\bV_{\bxi_i}^{1/2} \tbu)^\top(\bV_{\bxi_i}^{1/2} \tbu)
= \|\bV_{\bxi_i}^{1/2} \tbu\|^2$,
hence we obtain $\bV_{\bxi_i}^{1/2} \tbu= \bzero$, thus
$\bV_{\bxi_i} \tbu= \bV_{\bxi_i}^{1/2} (\bV_{\bxi_i}^{1/2} \tbu
) = \bzero$,
and hence $\tbu^\top\bV_{\bxi_i} = \bzero$.

First, we show \eqref{tZcond4}.
By \eqref{M2F},
\[
\sum_{k=1}^{{\lfloor nt\rfloor}} V_{k-1}^2
\tbu^\top\bV_{\bM_k} \tbu = \sum_{k=1}^{{\lfloor nt\rfloor}}
V_{k-1}^2 \tbu^\top\bV_\bvare \tbu,
\]
hence, in order to show \eqref{tZcond4}, it suffices to prove
\[
\sup_{t\in[0,T]} \Biggl| \frac{1}{n} \sum
_{k=1}^{{\lfloor nt\rfloor}} V_{k-1}^2 - t \EE
\bigl(\langle\tbu, \bvare_1 \rangle^2\bigr) \Biggr| \stoch0 .
\]
For all $k \in\NN$, by Remark~\ref{REMARK0},
$\langle\obV_\bxi\tbu, \tbu\rangle= 0$ implies
\begin{eqnarray*}
V_k &=& \sum_{j=1}^{X_{k-1,1}} (
\xi_{k,j,1,1} - \xi_{k,j,1,2}) + \sum_{j=1}^{X_{k-1,2}}
(\xi_{k,j,2,1} - \xi_{k,j,2,2}) + (\vare_{k,1} -
\vare_{k,2})
\\
&\ase&\vare_{k,1} - \vare_{k,2} = \langle\tbu,
\bvare_k \rangle.
\end{eqnarray*}
We have
\begin{eqnarray*}
\Biggl| \frac{1}{n} \sum_{k=1}^{{\lfloor nt\rfloor}}
V_{k-1}^2 - t \EE\bigl(\langle\tbu, \bvare_1
\rangle^2\bigr) \Biggr| &\leq&\frac{1}{n} \Biggl| \sum
_{k=1}^{{\lfloor nt\rfloor}} \bigl[ \langle\tbu, \bvare_{k-1}
\rangle^2 - \EE\bigl(\langle\tbu, \bvare_{k-1}
\rangle^2\bigr) \bigr] \Biggr|
\\
&&{} + \frac{|nt - {\lfloor nt\rfloor}|}{n} \EE\bigl(\langle\tbu, \bvare_k
\rangle^2\bigr) ,
\end{eqnarray*}
where $|nt - {\lfloor nt\rfloor}| \leq1$, hence, in order to show
\eqref{tZcond4}, it
suffices to prove
%
%e7.7 #&#
%e7.8 #&#
\begin{eqnarray}
\label{tZcond4m} %
&&\frac{1}{n} \sup_{t\in[0,T]} \Biggl|
\sum_{k=1}^{{\lfloor nt\rfloor}} \bigl[ \langle\tbu,
\bvare_k \rangle^2 - \EE\bigl(\langle\tbu,
\bvare_k \rangle^2\bigr) \bigr] \Biggr|
\nonumber\\[-8pt]\\[-8pt]
&&\quad = \frac{1}{n} \max_{N\in\{1,\ldots,{\lfloor nT\rfloor}\}} \Biggl| \sum
_{k=1}^N \bigl[ \langle\tbu, \bvare_k
\rangle^2 - \EE\bigl(\langle\tbu, \bvare_k
\rangle^2\bigr) \bigr] \Biggr| \stoch0 .\nonumber %
\end{eqnarray}
Applying Kolmogorov's maximal inequality, we obtain
\begin{eqnarray*}
&&\PP \Biggl( n^{-1} \max_{N\in\{1,\ldots,{\lfloor nT\rfloor}\}} \Biggl| \sum
_{k=1}^N \bigl[ \langle\tbu, \bvare_k
\rangle^2 - \EE\bigl(\langle\tbu, \bvare_k
\rangle^2\bigr) \bigr] \Biggr| \geq\vare \Biggr)
\\
&&\quad \leq\frac{1}{n^2 \vare^2} \var \Biggl( \sum_{k=1}^{\lfloor nT\rfloor}
\langle\tbu, \bvare_k \rangle^2 \Biggr) =
\frac{{\lfloor nT\rfloor}}{n^2 \vare^2} \var\bigl(\langle\tbu, \bvare_k
\rangle^2\bigr) \to0 ,\qquad  \mbox{as }n \to\infty
\end{eqnarray*}
for all $\vare> 0$, thus we conclude \eqref{tZcond4m}, and hence
\eqref{tZcond4}.

Now we turn to check \eqref{tZcond5}.
By \eqref{M2F},
\[
\sum_{k=1}^{{\lfloor nt\rfloor}} V_{k-1}
\tbu^\top\bV_{\bM_k} = \sum_{k=1}^{{\lfloor nt\rfloor}}
V_{k-1} \tbu^\top\bV_\bvare.
\]
Again by the strong law of large numbers,
$n^{-1} \sum_{k=1}^{{\lfloor nT\rfloor}} |V_{k-1}| \as t \EE
(|\langle\tbu, \bvare
_1 \rangle|)$
as $n \to\infty$ for all $T>0$, hence we conclude \eqref{tZcond5}.

Finally, we check \eqref{tZcond6}.
By \eqref{M2F},
\[
\sum_{k=1}^{{\lfloor nt\rfloor}} U_{k-1}
V_{k-1} \tbu^\top\EE\bigl( \bM _k
\bM_k^\top|\cF_{k-1} \bigr) = \sum
_{k=1}^{{\lfloor nt\rfloor}} U_{k-1} V_{k-1}
\tbu^\top\bV _\bvare.
\]
Applying $V_k = \langle\tbu, \bvare_k \rangle$, $k \in\NN$, and
Corollary~\ref{EEX_EEU_EEV}, we have
$\EE(|U_{k-1} V_{k-1}|) \leq\linebreak[4] \sqrt{\EE(U_{k-1}^2) \EE(V_{k-1}^2)} =
\OO(k)$,
which clearly implies \eqref{tZcond6}.
Condition (ii) of Theorem~\ref{Conv2DiffThm} can be checked again as
in case of
Theorem~\ref{main_conv}.
\end{pf*}

\begin{appendix}

%s8 #&#
\section{CLS estimators}
\label{section_estimators}

In order to analyse existence and uniqueness of the estimators given in
\eqref{CLSEr}, \eqref{CLSEd} and \eqref{CLSEab} in case of a
critical doubly
symmetric $2$-type Galton--Watson process, that is, when $\varrho= 1$,
we need
the following approximations.
%
%le8.1 #&#
\begin{Lem}\label{main_VV}
We have
\[
n^{-2} \Biggl( \sum_{k=1}^n
V_k^2 - \frac{\langle\obV_\bxi\tbu, \tbu\rangle}{4 \alpha\beta} \sum
_{k=1}^n U_{k-1} \Biggr) \stoch0 , \qquad \mbox{as }n \to\infty.
\]
\end{Lem}
\begin{pf}
In order to prove the statement, we derive a decomposition of
$\sum_{k=1}^n V_k^2$ as a sum of a martingale and some negligible terms.
Using recursion \eqref{rec_V}, Lemma~\ref{Moments} and \eqref{XUV},
we obtain
\begin{eqnarray*}
\EE\bigl(V_k^2 |\cF_{k-1}\bigr) &= &(\alpha-
\beta)^2 V_{k-1}^2 + 2 (\alpha- \beta) \langle
\tbu, {\mathbf{m}}_\bvare\rangle V_{k-1} + \langle\tbu, {
\mathbf{m}}_\bvare\rangle^2
\\
&&{} + \tbu^\top\EE\bigl(\bM_k \bM_k^\top|
\cF_{k-1}\bigr) \tbu
\\
&=& (\alpha- \beta)^2 V_{k-1}^2 +
\tfrac{1}{2} \tbu^\top(\bV_{\bxi_1} + \bV_{\bxi_2})
\tbu U_{k-1} + \mbox{constant} + \mbox{constant}\times V_{k-1}.
\end{eqnarray*}
Thus,
\begin{eqnarray*}
\sum_{k=1}^n V_k^2
&=& \sum_{k=1}^n \bigl[
V_k^2 - \EE\bigl(V_k^2 |
\cF_{k-1}\bigr) \bigr] + (\alpha- \beta)^2 \sum
_{k=1}^n V_{k-1}^2 +
\tbu^\top\obV_\bxi\tbu\sum_{k=1}^n
U_{k-1}
\\
&&{} + \OO(n) + \mbox{constant}\times\sum_{k=1}^n
V_{k-1}.
\end{eqnarray*}
Consequently,
%
%e8.1 #&#
\begin{eqnarray}
\label{sum_Vk2} %
 \sum_{k=1}^n
V_k^2 &=& \frac{1}{1 - (\alpha- \beta)^2} \sum
_{k=1}^n \bigl[ V_k^2 - \EE
\bigl(V_k^2 |\cF_{k-1}\bigr) \bigr]\nonumber
\\
&&{} + \frac{1}{1 - (\alpha- \beta)^2} \langle\obV_\bxi\tbu, \tbu\rangle\sum
_{k=1}^n U_{k-1}
\\
& &{}- \frac{(\alpha- \beta)^2}{1 - (\alpha- \beta)^2} V_n^2 + \OO(n) + \mbox{constant}
\times \sum_{k=1}^n V_{k-1}.\nonumber
\end{eqnarray}
Using \eqref{seged_UV_UNIFORM4} with $(\ell, i, j) = (8, 0, 2)$, we
obtain
\begin{eqnarray*}
\frac{1}{n^2} \sum_{k=1}^n \bigl[
V_k^2 - \EE\bigl(V_k^2 |
\cF_{k-1}\bigr) \bigr] \stoch0 , \qquad \mbox{as }n \to\infty.
\end{eqnarray*}
By Corollary~\ref{EEX_EEU_EEV}, we obtain $\EE(V_n^2) = \OO(n)$, and hence
$n^{-2} V_n^2 \stoch0$.
Moreover, $n^{-2}\times\allowbreak  \sum_{k=1}^n V_{k-1} \stoch0$ as $n \to\infty$
follows by \eqref{seged_UV_UNIFORM1} with the choices
$(\ell, i, j) = (4, 0, 1)$.
Consequently, by \eqref{sum_Vk2}, we obtain the statement, since
$1 - (\alpha- \beta)^2 = 4 \alpha\beta$.
\end{pf}
%
%le8.2 #&#
\begin{Lem}\label{main_VVt}
If $\langle\obV_\bxi\tbu, \tbu\rangle= 0$, then
\[
n^{-1} \sum_{k=1}^n
V_k^2 \as\EE\bigl(\langle\tbu, \bvare_1
\rangle^2\bigr) , \qquad \mbox{as }n \to\infty,
\]
and $\EE(\langle\tbu, \bvare_1 \rangle^2) = 0$ if and only if
$X_{k,1} \ase X_{k,2}$ for all $k \in\NN$.
\end{Lem}
\begin{pf}
By Remark~\ref{REMARK0}, $\langle\obV_\bxi\tbu, \tbu\rangle= 0$ implies
$V_k \ase\vare_{k,1} - \vare_{k,2} = \langle\tbu, \bvare_1
\rangle$ for
all $k \in\NN$, hence the convergence follows from the strong law of
large numbers.
Clearly $\EE(\langle\tbu, \bvare_1 \rangle^2) = 0$ is equivalent to
$\langle\tbu, \bvare_1 \rangle= \vare_{1,1} - \vare_{1,2} \ase
0$, and
hence it is equivalent to $X_{k,1} - X_{k,2} \ase0$ for all $k \in\NN$.
\end{pf}

Now we can prove existence and uniqueness of CLS estimators of the offspring
means and of the criticality parameter.
%
%pr8.3 #&#
\begin{Pro}\label{ExUn}
We have $\lim_{n \to\infty} \PP( (\bX_1, \ldots, \bX_n) \in H_n)
= 1$, where
$H_n$ is defined in \eqref{H_n}, and hence the probability of the
existence of a unique CLS estimator $\hvarrho_n$ converges to 1 as
$n \to\infty$, and this CLS estimator has the form given in \eqref{CLSEr}
whenever the sample $(\bX_1, \ldots, \bX_n)$ belongs to the set $H_n$.

If $\langle\obV_\bxi\tbu, \tbu\rangle> 0$, or if
$\langle\obV_\bxi\tbu, \tbu\rangle= 0$ and
$\EE(\langle\tbu, \bvare\rangle^2) > 0$, then
$\lim_{n \to\infty} \PP( (\bX_1, \ldots, \bX_n) \in\tH_n) =
1$, where
$\tH_n$ is defined in \eqref{tH_n}, and hence the probability of the
existence of unique CLS estimators $\hdelta_n$ and
$(\halpha_n, \hbeta_n)$ converges to 1 as $n \to\infty$.
The CLS estimator $\hdelta_n$ has the form given in \eqref{CLSEd}
whenever the sample $(\bX_1, \ldots, \bX_n)$ belongs to the set $\tH_n$.
The CLS estimator $(\halpha_n, \hbeta_n)$ has the form given in
\eqref{CLSEab} whenever the sample $(\bX_1, \ldots, \bX_n)$ belongs
to the
set $H_n \cap\tH_n$.
\end{Pro}
\begin{pf}
Recall convergence $\bcX^{(n)} \distr\bcX= \frac{1}{2} \cY\bone$ from
\eqref{convX}.
By Lemmas \ref{Conv2Funct} and \ref{Marci} one can show
%
%e8.2 #&#
\begin{eqnarray}
\label{seged2} \frac{1}{n^3} \sum_{k=1}^n
\bigl(X_{k-1,1}^2 + X_{k-1,2}^2\bigr)
\distr\frac{1}{2} \int_0^1
\cY_t^2 \,\mathrm{d}t ,\qquad  \mbox{as }n \to\infty,
\end{eqnarray}
see Isp\'any \textit{et al.} \cite{IKP}, Proposition A.4.
Since ${\mathbf{m}}_\bvare\ne\bzero$, by the SDE \eqref{Y}, we have
$\PP ( \cY_t=0, t \in[0,1]  ) = 0$, which implies that
$\PP ( \int_0^1 \cY_t^2 \,\mathrm{d}t > 0  ) = 1$.
Consequently, the distribution function of $\int_0^1 \cY_t^2 \,\mathrm
{d}t$ is
continuous at 0, and hence, by \eqref{seged2},
\[
\PP \Biggl( \sum_{k=1}^n \langle\bone,
\bX_{k-1}\rangle^2 > 0 \Biggr) \to\PP \biggl(
\frac{1}{2} \int_0^1
\cY_t^2 \,\mathrm{d}t > 0 \biggr) = 1 ,\qquad  \mbox{as }n \to
\infty.
\]
Now suppose that $\langle\obV_\bxi\tbu, \tbu\rangle> 0$ holds.
In a similar way, using Lemma~\ref{main_VV}, convergence \eqref
{convX}, and
Lemmas \ref{Conv2Funct} and \ref{Marci}, one can show
\[
\frac{1}{n^2} \sum_{k=1}^n \langle
\tbu, \bX_{k-1}\rangle^2 \distr \frac{\langle\obV_\bxi\tbu, \tbu\rangle}{4 \alpha\beta} \int
_0^1 \cY_t \,\mathrm{d}t , \qquad \mbox{as
}n \to\infty,
\]
implying
\[
\PP \Biggl( \sum_{k=1}^n \langle\tbu,
\bX_{k-1}\rangle^2 > 0 \Biggr) \to\PP \biggl( \int
_0^1 \cY_t \,\mathrm{d}t > 0 \biggr)
= 1 , \qquad \mbox{as }n \to\infty,
\]
hence we obtain the statement under the assumption
$\langle\obV_\bxi\tbu, \tbu\rangle> 0$.

Next, we suppose that $\langle\obV_\bxi\tbu, \tbu\rangle= 0$ and
$\EE(\langle\tbu, \bvare\rangle^2) > 0$ hold.
Then
\[
\PP \Biggl( \sum_{k=1}^n \langle\tbu,
\bX_{k-1}\rangle^2 > 0 \Biggr) = \PP \Biggl(
\frac{1}{n} \sum_{k=1}^n
V_{k-1}^2 > 0 \Biggr) \to1 ,\qquad  \mbox{as }n \to\infty,
\]
since Lemma~\ref{main_VVt} yields
$n^{-1} \sum_{k=1}^n V_{k-1}^2
\stoch\EE(\langle\tbu, \bvare_1 \rangle^2) > 0$,
and hence we conclude the statement under the assumptions
$\langle\obV_\bxi\tbu, \tbu\rangle= 0$ and
$\EE(\langle\tbu, \bvare\rangle^2) > 0$.
\end{pf}

%s9 #&#
\section{Estimations of moments}
\label{section_moments}

In the proof of Theorem~\ref{main}, good bounds for moments of the random
vectors and variables $(\bM_k)_{k\in\ZZ_+}$, $(\bX_k)_{k\in\ZZ_+}$,
$(U_k)_{k\in\ZZ_+}$ and $(V_k)_{k\in\ZZ_+}$ are extensively used.
First note that, for all $k \in\NN$, $\EE( \bM_k |\cF_{k-1} )
= \bzero$
and $\EE(\bM_k) = \bzero$, since
$\bM_k = \bX_k - \EE(\bX_k |\cF_{k-1})$.
%
%le9.1 #&#
\begin{Lem}\label{Moments}
Let $(\bX_k)_{k\in\ZZ_+}$ be a $2$-type Galton--Watson process with immigration
and with\linebreak[4]  $\bX_0 = \bzero$.
If $\EE(\|\bxi_{1,1,1}\|^2) < \infty$, $\EE(\|\bxi_{1,1,2}\|^2) <
\infty$
and $\EE(\|\bvare_1\|^2) < \infty$ then
%
%e9.1 #&#
\begin{equation}
\label{Mcond} \EE\bigl( \bM_k \bM_k^\top|
\cF_{k-1} \bigr) = X_{k-1,1} \bV_{\bxi_1} +
X_{k-1,2} \bV_{\bxi_2} + \bV_\bvare,\qquad  k \in\NN.
\end{equation}
If $\EE(\|\bxi_{1,1,1}\|^3) < \infty$, $\EE(\|\bxi_{1,1,2}\|^3) <
\infty$
and $\EE(\|\bvare_1\|^3) < \infty$, then
%
%e9.2 #&#
%e9.3 #&#
\begin{eqnarray}
\label{M3cond} %
\EE\bigl( \bM_k^{\otimes3} |
\cF_{k-1} \bigr) &=& X_{k-1,1} \EE\bigl[(\bxi_{1,1,1} -
\EE(\bxi_{1,1,1})^{\otimes3}\bigr]
\nonumber\\[-8pt]\\[-8pt]
&& {} + X_{k-1,2} \EE\bigl[(\bxi_{1,1,2} - \EE(\bxi_{1,1,2})^{\otimes3}
\bigr] + \EE\bigl[(\bvare_1 - \EE(\bvare_1)^{\otimes3}
\bigr] ,\qquad  k \in\NN.\nonumber %
\end{eqnarray}
\end{Lem}
\begin{pf} By \eqref{GWI(2)} and \eqref{Mk}, $\bM_k$ has the
form
%
%e9.4 #&#
\begin{equation}
\label{Mdeco} \sum_{j=1}^{X_{k-1,1}} \bigl(
\bxi_{k,j,1} - \EE(\bxi_{k,j,1}) \bigr) + \sum
_{j=1}^{X_{k-1,2}} \bigl( \bxi_{k,j,2} - \EE(
\bxi_{k,j,2}) \bigr) + \bigl( \bvare_k - \EE(
\bvare_k) \bigr)
\end{equation}
for all $k \in\NN$.
The random vectors
$ \{\bxi_{k,j,1} - \EE(\bxi_{k,j,1}), \bxi_{k,j.2} - \EE(\bxi
_{k,j,2}),
\bvare_k - \EE(\bvare_k)
: j \in\NN \}$
are independent of each other, independent of $\cF_{k-1}$, and have
zero mean vector, thus we conclude \eqref{Mcond} and \eqref{M3cond}.
\end{pf}
%
%le9.2 #&#
\begin{Lem}\label{LEM_moments_seged0}
Let $(\bzeta_k)_{k\in\NN}$ be independent and identically distributed
random vectors with values in $\RR^d$ such that
$\EE(\|\bzeta_1\|^\ell) < \infty$ with some $\ell\in\NN$.
\begin{enumerate}[(ii)]
\item[(i)]
Then there exists $\bQ= (Q_1, \ldots, Q_{d^\ell})\dvtx  \RR\to\RR
^{d^\ell}$,
where $Q_1, \ldots, Q_{d^\ell}$ are polynomials having degree at most
$\ell-1$ such that
\[
\EE \bigl((\bzeta_1 + \cdots+ \bzeta_N)^{\otimes\ell}
\bigr) = N^\ell \bigl[\EE(\bzeta_1) \bigr]^{\otimes\ell} +
\bQ(N) ,\qquad  N \in\NN, N \geq\ell.
\]
\item[(ii)]
If $\EE(\bzeta_1) = \bzero$, then there exists
$\bR= (R_1, \ldots, R_{d^\ell}) \dvtx  \RR\to\RR^{d^\ell}$, where
$R_1$, \ldots, $R_{d^\ell}$ are polynomials having degree at most
$\lfloor\ell/2 \rfloor$ such that
\[
\EE \bigl((\bzeta_1 + \cdots+ \bzeta_N)^{\otimes\ell}
\bigr) = \bR(N) ,\qquad  N \in\NN, N \geq\ell.
\]
\end{enumerate}
The coefficients of the polynomials $\bQ$ and $\bR$ depend on the
moments $\EE(\bzeta_{i_1} \otimes\cdots\otimes\bzeta_{i_\ell})$,
$i_1, \ldots, i_\ell\in\{1, \ldots, N\}$.
\end{Lem}
\begin{pf}
(i) We have
\begin{eqnarray*}
&&\EE \bigl((\bzeta_1 + \cdots+ \bzeta_N)^{\otimes\ell}
\bigr)
\\
&&\quad = \mathop{\mathop{\sum}_{s \in\{1,\ldots,\ell\}, k_1, \ldots,
k_s \in\ZZ_+,}}_{ k_1 + 2 k_2 + \cdots+ s k_s = \ell, k_s \ne0}
\binom{N} {k_1} \binom{N - k_1} {k_2} \cdots
\binom{N - k_1 - \cdots- k_{s-1}} {k_s}
\\
&&\qquad {}\times\sum_{(i_1,\ldots,i_\ell) \in P_{k_1,\ldots,k_s}^{(N,\ell)}} \EE(\bzeta_{i_1}
\otimes\cdots\otimes\bzeta_{i_\ell}) , %
\end{eqnarray*}
where the set $P_{k_1,\ldots,k_s}^{(N,\ell)}$ consists of
permutations of all
the multisets containing pairwise different elements
$j_{k_1}, \ldots, j_{k_s}$ of the set $\{1, \ldots, N\}$ with
multiplicities $k_1, \ldots, k_s$, respectively.
Since
\begin{eqnarray*}
&&\binom{N} {k_1}\binom{N - k_1} {k_2} \cdots
\binom{N - k_1 - \cdots- k_{s-1}} {k_s}
\\
&&\quad = \frac{N (N-1) \cdots(N - k_1 - k_2 - \cdots- k_s + 1)} {
k_1! k_2! \cdots k_s!}
\end{eqnarray*}
is a polynomial of the variable $N$ having degree
$k_1 + \cdots+ k_s \leq\ell$, there exists
$\bP= (P_1, \ldots, P_{d^\ell}) \dvtx  \RR\to\RR^{d^\ell}$, where
$P_1, \ldots, P_{d^\ell}$ are polynomials having degree at most
$\ell$ such that
$\EE ((\bzeta_1 + \cdots+ \bzeta_N)^{\otimes\ell} ) =
\bP(N)$.
A term of degree $\ell$ can occur only in case
$k_1 + \cdots+ k_s = \ell$, when $k_1 + 2k_2 + \cdots+ sk_s = \ell$
implies $s = 1$ and $k_1 = \ell$, thus the corresponding term of
degree $\ell$ is
$N (N-1) \cdots(N-\ell+1)  [\EE(\bzeta_1) ]^{\otimes\ell
}$, hence
we obtain the statement.
Part (ii) can be proved in a similar way.
\end{pf}

Lemma~\ref{LEM_moments_seged0} can be generalized in the following way.
%
%le9.3 #&#
\begin{Lem}\label{LEM_moments_seged}
For each $i \in\NN$, let $(\bzeta_{i,k})_{k\in\NN}$ be independent and
identically distributed random vectors with values in $\RR^d$ such that
$\EE(\|\bzeta_{i,1}\|^\ell) < \infty$ with some $\ell\in\NN$.
Let $j_1, \ldots, j_\ell\in\NN$.
\begin{enumerate}[(ii)]
\item[(i)]
Then there exists
$\bQ= (Q_1, \ldots, Q_{d^\ell}) \dvtx  \RR^\ell\to\RR^{d^\ell}$,
where $Q_1, \ldots, Q_{d^\ell}$ are polynomials of $\ell$
variables having degree at most $\ell-1$ such that
\begin{eqnarray*}
&&\EE \bigl((\bzeta_{j_1,1} + \cdots+ \bzeta_{j_1,N_1}) \otimes
\cdots\otimes (\bzeta_{j_\ell,1} + \cdots+ \bzeta_{j_\ell,N_\ell}) \bigr)
\\
&&\quad = N_1 \cdots N_\ell \EE(\bzeta_{j_1,1}) \otimes
\cdots\otimes\EE(\bzeta_{j_\ell,1}) + \bQ(N_1, \ldots,
N_\ell)
\end{eqnarray*}
for $N_1, \ldots, N_\ell\in\NN$ with $N_1 \geq\ell$, \ldots,
$N_\ell\geq\ell$.
\item[(ii)]
If $\EE(\bzeta_{j_1,1}) = \cdots= \EE(\bzeta_{j_\ell,1}) = \bzero
$, then
there exists $\bR= (R_1, \ldots, R_{d^\ell}) \dvtx  \RR^\ell\to\RR
^{d^\ell}$,
where\vadjust{\goodbreak} $R_1$, \ldots, $R_{d^\ell}$ are polynomials of $\ell$
variables having degree at most $\lfloor\ell/2 \rfloor$ such that
\[
\EE \bigl((\bzeta_{j_1,1} + \cdots+ \bzeta_{j_1,N_1}) \otimes
\cdots\otimes (\bzeta_{j_\ell,1} + \cdots+ \bzeta_{j_\ell,N_\ell}) \bigr)
\\
= \bR(N_1, \ldots, N_\ell)
\]
for $N_1, \ldots, N_\ell\in\NN$ with $N_1 \geq\ell$, \ldots,
$N_\ell\geq\ell$.
\end{enumerate}
The coefficients of the polynomials $\bQ$ and $\bR$ depend on the
moments $\EE(\bzeta_{j_1,i_1} \otimes\cdots\otimes\bzeta_{j_\ell
,i_\ell})$,
$i_1 \in\{1, \ldots, N_1\}$, \ldots, $i_\ell\in\{1, \ldots,
N_\ell\}$.
\end{Lem}
%
%le9.4 #&#
\begin{Lem}\label{LEM_Putzer}
If $(\alpha, \beta) \in[0, 1]$ with $\alpha+ \beta= 1$, then the
matrix ${\mathbf{m}}_\bxi$ defined in \eqref{bA} has eigenvalues
$1$ and
$\alpha- \beta$, and the powers of ${\mathbf{m}}_\bxi$ take the form
\begin{eqnarray*}
{\mathbf{m}}_\bxi^j = \frac{1}{2} %
\lleft[\matrix{ 1 & 1
\cr
1 & 1 } \rright] %
+ \frac{1}{2} (
\alpha- \beta)^j %
\lleft[\matrix{ 1 & -1
\cr
-1 & 1 }
\rright] %
,\qquad  j \in\ZZ_+ .
\end{eqnarray*}
Consequently, $\|{\mathbf{m}}_\bxi^j\| = \OO(1)$, that is,
$\sup_{j \in\NN} \|{\mathbf{m}}_\bxi^j\| < \infty$.
\end{Lem}
%
%le9.5 #&#
\begin{Lem}\label{LEM_moments_X}
Let $(\bX_k)_{k\in\ZZ_+}$ be a $2$-type doubly symmetric Galton--Watson
process with immigration with offspring means $(\alpha, \beta) \in
[0, 1]$
such that $\alpha+ \beta= 1$ (hence it is critical).
Suppose $\bX_0 = \bzero$, and $\EE(\|\bxi_{1,1,1}\|^\ell) < \infty$,
$\EE(\|\bxi_{1,1,2}\|^\ell) < \infty$, $\EE(\|\bvare_1\|^\ell) <
\infty$
with some $\ell\in\NN$.
Then $\EE(\|\bX_k\|^\ell) = \OO(k^\ell)$, that is,
$\sup_{k \in\NN} k^{-\ell} \EE(\|\bX_k\|^\ell) < \infty$.
\end{Lem}
\begin{pf} The statement is clearly equivalent with
$\EE (|P(X_{k,1},X_{k,2})| ) \leq c_P k^\ell$, $k \in\NN
$, for
all polynomials $P$ of two variables having degree at most $\ell$,
where $c_P$ depends only on $P$.

If $\ell= 1$, then \eqref{EXk} and Lemma~\ref{LEM_Putzer} imply
\begin{eqnarray*}
\EE(\bX_k) = \sum_{j=0}^{k-1} {
\mathbf{m}}_\bxi^j {\mathbf{m}}_\bvare = \lleft(\frac{k}{2} %
\lleft[\matrix{ 1 & 1
\cr
1 & 1 } \rright]
+ \frac{1 - (\alpha-\beta)^k}{4\beta} %
\lleft[\matrix{ 1 & -1
\cr
-1 & 1 }
\rright] %
\rright) {\mathbf{m}}_\bvare,
\end{eqnarray*}
for all $k \in\NN$, which yields the statement.

Using part (i) of Lemma~\ref{LEM_moments_seged} and separating the
terms having
degree 2 and less than 2, we obtain
\begin{eqnarray*}
\EE\bigl(\bX_k^{\otimes2} |\cF_{k-1}\bigr)
&=& X_{k-1,1}^2 {\mathbf{m}}_{\bxi_1}^{\otimes2} +
X_{k-1,2}^2 {\mathbf{m}}_{\bxi_2}^{\otimes2} +
X_{k-1,1} X_{k-1,2} ({\mathbf{m}}_{\bxi_1} \otimes{
\mathbf{m}}_{\bxi_2} + {\mathbf{m}}_{\bxi_2} \otimes{
\mathbf{m}}_{\bxi_1})
\\
&& {} + \bQ_2(X_{k-1,1},X_{k-1,2})
\\
& =& (X_{k-1,1} {\mathbf{m}}_{\bxi_1} + X_{k-1,2} {
\mathbf {m}}_{\bxi_2})^{\otimes2} + \bQ_2(X_{k-1,1},X_{k-1,2})
\\
& =& ({\mathbf{m}}_\bxi\bX_{k-1})^{\otimes2} +
\bQ_2(X_{k-1,1},X_{k-1,2}) = {\mathbf{m}}_\bxi^{\otimes2}
\bX_{k-1}^{\otimes2} + \bQ _2(X_{k-1,1},X_{k-1,2})
,
\end{eqnarray*}
where $\bQ_2 = (Q_{2,1}, Q_{2,2}, Q_{2,3}, Q_{2,4}) \dvtx  \RR^2 \to\RR
^4$, and
$Q_{2,1}$, $Q_{2,2}$, $Q_{2,3}$ and $Q_{2,4}$ are polynomials of two
variables having degree at most $1$.
Hence
\[
\EE\bigl(\bX_k^{\otimes2}\bigr) = {\mathbf{m}}_\bxi^{\otimes2}
\EE\bigl(\bX_{k-1}^{\otimes2}\bigr) + \EE\bigl[\bQ_2(X_{k-1,1},X_{k-1,2})
\bigr] .
\]
In a similar way,
\[
\EE\bigl(\bX_k^{\otimes\ell}\bigr) = {\mathbf{m}}_\bxi^{\otimes\ell}
\EE\bigl(\bX_{k-1}^{\otimes\ell}\bigr) + \EE\bigl[\bQ_\ell(X_{k-1,1},X_{k-1,2})
\bigr] ,
\]
where $\bQ_\ell= (Q_{\ell,1}, \ldots, Q_{\ell,2^\ell}) \dvtx  \RR^2
\to\RR^{2^\ell}$,
and $Q_{\ell,1}, \ldots, Q_{\ell,2^\ell}$ are polynomials of two variables
having degree at most $\ell-1$, implying
\begin{eqnarray*}
\EE\bigl(\bX_k^{\otimes\ell}\bigr) &=& \sum
_{j=1}^k \bigl({\mathbf{m}}_\bxi^{\otimes\ell}
\bigr)^{k-j} \EE \bigl[\bQ_\ell(X_{j-1,1},X_{j-1,2})
\bigr]
\\
&=& \sum_{j=0}^{k-1} \bigl({
\mathbf{m}}_\bxi^{\otimes\ell}\bigr)^j \EE\bigl[
\bQ_\ell (X_{k-j-1,1},X_{k-j-1,2})\bigr]
\\
&= &\sum_{j=0}^{k-1} \bigl({
\mathbf{m}}_\bxi^j\bigr)^{\otimes\ell} \EE\bigl[
\bQ_\ell (X_{k-j-1,1},X_{k-j-1,2})\bigr] .
\end{eqnarray*}
Let us suppose now that the statement holds for $1, \ldots, \ell-1$.
Then
\[
\EE\bigl[\bigl|Q_{\ell,i}(X_{k-j-1,1},X_{k-j-1,2})\bigr|\bigr] \leq
c_{Q_{\ell,i}} k^{\ell-1} ,\qquad  k \in\NN,\ i \in\bigl\{1, \ldots,
2^\ell\bigr\} .
\]
By Lemma~\ref{LEM_Putzer} $\|({\mathbf{m}}_\bxi^j)^{\otimes\ell
}\| = \OO(1)$, hence we
obtain the assertion for $\ell$.
\end{pf}
%
%co9.6 #&#
\begin{Cor}\label{EEX_EEU_EEV}
Let $(\bX_k)_{k\in\ZZ_+}$ be a $2$-type doubly symmetric Galton--Watson
process with immigration having offspring means
$(\alpha, \beta) \in(0, 1)^2$
such that $\alpha+ \beta= 1$ (hence it is critical and positively
regular).
Suppose $\bX_0 = \bzero$, and $\EE(\|\bxi_{1,1,1}\|^\ell) < \infty$,
$\EE(\|\bxi_{1,1,2}\|^\ell) < \infty$, $\EE(\|\bvare_1\|^\ell) <
\infty$
with some $\ell\in\NN$.
Then $\EE(\|\bX_k\|^\ell) = \OO(k^\ell)$,
$\EE(\bM_k^{\otimes\ell}) = \OO(k^{\lfloor\ell/2 \rfloor})$,
$\EE(U^\ell_k ) = \OO(k^\ell)$ and $\EE(V^{2j}_k ) = \OO(k^j)$ for
$j \in\ZZ_+$ with $2 j \leq\ell$.
\end{Cor}
\begin{pf}
The first statement is just Lemma~\ref{LEM_moments_X}.
Next, we turn to prove $\EE(\bM_k^{\otimes\ell}) = \OO(k^{\lfloor
\ell/2 \rfloor})$.
Using \eqref{Mdeco}, part (ii) of Lemma~\ref{LEM_moments_seged}, and
that the
random vectors
$ \{\bxi_{k,j,1} - \EE(\bxi_{k,j,1}), \bxi_{k,j.2} - \EE(\bxi
_{k,j,2}),
\bvare_k - \EE(\bvare_k)
: j \in\NN \}$
are independent of each other, independent of $\cF_{k-1}$, and have
zero mean vector, we obtain
$\EE(\bM_k^{\otimes\ell} |\cF_{k-1}) = \bR(X_{k-1,1},X_{k-1,2})$
with $\bR= (R_1, \ldots, R_{2^\ell}) \dvtx  \RR^2 \to\RR^{2 \ell}$, where
$R_1, \ldots, R_{2^\ell}$ are polynomials of two variables having degree
at most $\ell/2$.
Hence $\EE(\bM_k^{\otimes\ell}) = \EE(\bR(X_{k-1,1},X_{k-1,2}))$.
By Lemma~\ref{LEM_moments_X}, we conclude
$\EE(\bM_k^{\otimes\ell}) = \OO(k^{\lfloor\ell/2 \rfloor})$.
The rest of the proof can be carried out as in Corollary~9.1 of
Barczy \textit{et al.} \cite{BarIspPap2}.
\end{pf}

The next corollary can be derived as Corollary~9.2 of
Barczy \textit{et al.} \cite{BarIspPap2}.
%
%co9.7 #&#
\begin{Cor}\label{LEM_UV_UNIFORM}
Let $(\bX_k)_{k\in\ZZ_+}$ be a $2$-type doubly symmetric Galton--Watson
process with immigration having offspring means
$(\alpha, \beta) \in(0, 1)^2$
such that $\alpha+ \beta= 1$ (hence, it is critical and positively
regular).
Suppose $\bX_0 = \bzero$, and $\EE(\|\bxi_{1,1,1}\|^\ell) < \infty$,
$\EE(\|\bxi_{1,1,2}\|^\ell) < \infty$, $\EE(\|\bvare_1\|^\ell) <
\infty$
with some $\ell\in\NN$.
Then
\begin{enumerate}[(iii)]
\item[(i)]
for all $i,j\in\ZZ_+$ with $\max\{i,j\} \leq\lfloor\ell/2 \rfloor$,
and for all $\kappa> i + \frac{j}{2} + 1$, we have
%
%e9.5 #&#
\begin{eqnarray}
\label{seged_UV_UNIFORM1} n^{-\kappa} \sum_{k=1}^n
\bigl\vert U_k^i V_k^j\bigr\vert \stoch0 ,\qquad
\mbox{as }n\to\infty,
\end{eqnarray}
\item[(ii)]
for all $i,j\in\ZZ_+$ with $\max\{i,j\} \leq\ell$, for all $T>0$,
and for all $\kappa> i + \frac{j}{2} + \frac{i+j}{\ell}$, we have
%
%e9.6 #&#
\begin{eqnarray}
\label{seged_UV_UNIFORM2} n^{-\kappa} \sup_{t\in[0,T]} \bigl\vert
U_{\lfloor nt\rfloor}^i V_{\lfloor nt\rfloor}^j \bigr\vert\stoch0 ,\qquad
\mbox{as }n\to\infty,
\end{eqnarray}
\item[(iii)]
for all $i,j\in\ZZ_+$ with $\max\{i,j\} \leq\lfloor\ell/4 \rfloor$,
for all $T>0$, and for all
$\kappa> i + \frac{j}{2} + \frac{1}{2}$, we have
%
%e9.7 #&#
\begin{eqnarray}
\label{seged_UV_UNIFORM4} n^{-\kappa} \sup_{t\in[0,T]} \Biggl\llvert \sum
_{k=1}^{\lfloor nt\rfloor}\bigl[U_k^i
V_k^j - \EE\bigl(U_k^i
V_k^j |\cF_{k-1}\bigr)\bigr] \Biggr\rrvert
\stoch0 , \qquad \mbox{as }n\to\infty.
\end{eqnarray}
\end{enumerate}
\end{Cor}
%
%re9.8 #&#
\begin{Rem}
In the special case $(\ell, i, j) = (2, 1, 0)$, one can improve
\eqref{seged_UV_UNIFORM2}, namely, one can show
%
%e9.8 #&#
\begin{equation}
\label{seged_UV_UNIFORM5} n^{-\kappa} \sup_{t\in[0,T]} U_{{\lfloor nt\rfloor}}
\stoch0 ,\qquad  \mbox{as }n\to\infty\mbox{ for }\kappa> 1,
\end{equation}
see Barczy \textit{et al.} \cite{BarIspPap2}.
\end{Rem}

%s10 #&#
\section{A version of the continuous mapping theorem}
\label{app_C}

A function $f \dvtx  \RR_+ \to\RR^d$ is called \emph{c\`adl\`ag} if it
is right
continuous with left limits.
Let $\DD(\RR_+, \RR^d)$ and $\CC(\RR_+, \RR^d)$ denote the space of
all $\RR^d$-valued c\`adl\`ag and continuous functions on $\RR_+$,
respectively.
Let $\cB(\DD(\RR_+, \RR^d))$ denote the Borel $\sigma$-algebra on
$\DD(\RR_+, \RR^d)$ for the metric defined in Jacod and Shiryaev
\cite{JSh}, Chapter VI, (1.26) (with this metric $\DD(\RR_+, \RR
^d)$ is a
complete and separable metric space and the topology induced by this
metric is
the so-called Skorokhod topology).
For $\RR^d$-valued stochastic processes $(\bcY_t)_{t \in\RR_+}$ and
$(\bcY^{(n)}_t)_{t \in\RR_+}$, $n \in\NN$, with c\`adl\`ag paths,
we write
$\bcY^{(n)} \distr\bcY$ if the distribution of $\bcY^{(n)}$ on the
space $(\DD(\RR_+, \RR), \cB(\DD(\RR_+, \RR^d)))$ converges
weakly to the
distribution of $\bcY$ on the space
$(\DD(\RR_+, \RR), \cB(\DD(\RR_+, \RR^d)))$ as $n \to\infty$.
Concerning the notation $\distr$ we note that if $\xi$ and $\xi_n$,
$n \in\NN$, are random elements with values in a metric space $(E, d)$,
then we also denote by $\xi_n \distr\xi$ the weak convergence of the
distributions of $\xi_n$ on the space $(E, \cB(E))$ towards the
distribution of $\xi$ on the space $(E, \cB(E))$ as $n \to\infty$,
where $\cB(E)$ denotes the Borel $\sigma$-algebra on $E$ induced by
the given metric $d$.

The following version of continuous mapping theorem can be found, for example,
in Kallenberg \cite[Theorem~3.27]{K}.
%
%le10.1 #&#
\begin{Lem}\label{Lem_Kallenberg}
Let $(S, d_S)$ and $(T, d_T)$ be metric spaces and
$(\xi_n)_{n \in\NN}$, $\xi$ be random elements with values in $S$
such that $\xi_n \distr\xi$ as $n \to\infty$.
Let $f \dvtx  S \to T$ and $f_n \dvtx  S \to T$, $n \in\NN$, be measurable
mappings and $C \in\cB(S)$ such that $\PP(\xi\in C) = 1$ and
$\lim_{n \to\infty} d_T(f_n(s_n), f(s)) = 0$ if
$\lim_{n \to\infty} d_S(s_n,s) = 0$ and $s \in C$.
Then $f_n(\xi_n) \distr f(\xi)$, as $n \to\infty$.
\end{Lem}

For the case $S = \DD(\RR_+, \RR^d)$ and $T = \RR^q$
(or $T = \DD(\RR_+,\RR^q)$), where $d$, $q \in\NN$, we formulate
a consequence of Lemma~\ref{Lem_Kallenberg}.

For functions $f$ and $f_n$, $n \in\NN$, in $\DD(\RR_+, \RR^d)$,
we write $f_n \lu f$ if $(f_n)_{n \in\NN}$ converges to $f$
locally uniformly, that is, if $\sup_{t \in[0,T]} \|f_n(t) - f(t)\| \to
0$ as
$n \to\infty$ for all $T > 0$.
For measurable mappings
$\Phi\dvtx  \DD(\RR_+, \RR^d) \to\RR^q$
(or $\Phi\dvtx  \DD(\RR_+, \RR^d) \to\DD(\RR_+,\RR^q)$) and
$\Phi_n \dvtx  \DD(\RR_+, \RR^d) \to\RR^q$
(or $\Phi_n \dvtx  \DD(\RR_+, \RR^d) \to\DD(\RR_+,\RR^q)$), $n \in
\NN$,
we will denote by $C_{\Phi, (\Phi_n)_{n \in\NN}}$ the set of all functions
$f \in\CC(\RR_+, \RR^d)$ such that
$\Phi_n(f_n) \to\Phi(f)$ (or $\Phi_n(f_n) \to\lu\Phi(f)$) whenever
$f_n \lu f$ with $f_n \in\DD(\RR_+, \RR^d)$, $n \in\NN$.

We will use the following version of the continuous mapping theorem several
times, see, for example, Isp\'any and Pap \cite{IspPap}, Lemma~3.1.
%
%le10.2 #&#
\begin{Lem}\label{Conv2Funct}
Let $d,q\in\NN$, and $(\bcU_t)_{t \in\RR_+}$ and
$(\bcU^{(n)}_t)_{t \in\RR_+}$, $n \in\NN$, be $\RR^d$-valued stochastic
processes with c\`adl\`ag paths such that $\bcU^{(n)} \distr\bcU$.
Let $\Phi\dvtx  \DD(\RR_+, \RR^d) \to\RR^q$
(or $\Phi\dvtx  \DD(\RR_+, \RR^d) \to\DD(\RR_+,\RR^q)$) and
$\Phi_n \dvtx  \DD(\RR_+, \RR^d) \to\RR^q$
(or $\Phi_n \dvtx  \DD(\RR_+, \RR^d) \to\DD(\RR_+,\RR^q)$), $n \in
\NN$,
be measurable mappings such that there exists
$C \subset C_{\Phi,(\Phi_n)_{n\in\NN}}$ with $C \in\cB(\DD(\RR
_+, \RR^d))$
and $\PP(\bcU\in C) = 1$.
Then $\Phi_n(\bcU^{(n)}) \distr\Phi(\bcU)$.
\end{Lem}

In order to apply Lemma~\ref{Conv2Funct}, we will use the following statement
several times, see Barczy \textit{et al.} \cite{BarIspPap2}, Lemma B.3.
%
%le10.3 #&#
\begin{Lem}\label{Marci}
Let $d, p, q \in\NN$, $h\dvtx \RR^d\to\RR^q$ be a continuous function and
$K \dvtx  [0,1] \times\RR^{2d} \to\RR^p$ be a function such that for all
$R > 0$ there exists $C_R > 0$ such that
%
%e10.1 #&#
\begin{equation}
\label{Lipschitz} \bigl\| K(s, x) - K(t, y) \bigr\| \leq C_R \bigl( | t - s | + \| x
- y \| \bigr)
\end{equation}
for all $s, t \in[0, 1]$ and $x, y \in\RR^{2d}$ with
$\| x \| \leq R$ and $\| y \| \leq R$.
Moreover, let us define the mappings
$\Phi, \Phi_n \dvtx  \DD(\RR_+, \RR^d) \to\RR^{q+p}$, $n \in\NN$, by
\begin{eqnarray*}
\Phi_n(f) &:=& \Biggl( h\bigl(f(1)\bigr), \frac{1}{n} \sum
_{k=1}^n K \biggl( \frac{k}{n}, f
\biggl( \frac{k}{n} \biggr), f \biggl( \frac{k-1}{n} \biggr) \biggr)
\Biggr) ,
\\
\Phi(f) & := &\biggl( h\bigl(f(1)\bigr), \int_0^1
K\bigl( u, f(u), f(u) \bigr) \,\mathrm{d}u \biggr)
\end{eqnarray*}
for all $f \in\DD(\RR_+, \RR^d)$.
Then the mappings $\Phi$ and $\Phi_n$, $n \in\NN$, are measurable,
and $C_{\Phi,(\Phi_n)_{n \in\NN}} = \CC(\RR_+, \RR^d) \in\cB
(\DD(\RR_+, \RR^d))$.
\end{Lem}

%s11 #&#
\section{Convergence of random step processes}
\label{section_conv_step_processes}

We recall a result about convergence of random step processes towards a
diffusion process, see Isp\'any and Pap \cite{IspPap}.
This result is used for the proof of convergence \eqref{conv_Z}.
%
%th11.1 #&#
\begin{Thm}\label{Conv2DiffThm}
Let $\bgamma\dvtx  \RR_+ \times\RR^d \to\RR^{d \times r}$ be a continuous
function.
Assume that uniqueness in the sense of probability law holds for the SDE
%
%e11.1 #&#
\begin{equation}
\label{SDE} \mathrm{d} \bcU_t = \bgamma(t, \bcU_t)
\,\mathrm{d}\bcW_t , \qquad t \in\RR_+,
\end{equation}
with initial value $\bcU_0 = {\mathbf{u}}_0$ for all
${\mathbf{u}}_0 \in\RR^d$, where
$(\bcW_t)_{t \in\RR_+}$ is an $r$-dimensional standard Wiener process.
Let $(\bcU_t)_{t \in\RR_+}$ be a solution of \textup{\eqref{SDE}} with
initial value
$\bcU_0 = \bzero\in\RR^d$.

For each $n \in\NN$, let $(\bU^{(n)}_k)_{k \in\NN}$ be a sequence of
$d$-dimensional martingale differences with respect to a filtration
$(\cF^{(n)}_k)_{k \in\ZZ_+}$, that is, $\EE(\bU^{(n)}_k |\cF
^{(n)}_{k-1}) = 0$,
$n \in\NN$, $k \in\NN$.
Let
\[
\bcU^{(n)}_t := \sum_{k=1}^{{\lfloor nt\rfloor}}
\bU^{(n)}_k ,\qquad  t \in\RR_+,\ n \in\NN.
\]
Suppose $\EE ( \|\bU^{(n)}_k\|^2  ) < \infty$ for all
$n, k \in\NN$.
Suppose that for each $T > 0$,
\begin{enumerate}[(ii)]
\item[(i)]
$\sup_{t\in[0,T]}
\llVert  \sum_{k=1}^{{\lfloor nt\rfloor}}
\EE (\bU^{(n)}_k (\bU^{(n)}_k)^\top|\cF^{(n)}_{k-1} )
- \int_0^t
\bgamma(s,\bcU^{(n)}_s) \bgamma(s,\bcU^{(n)}_s)^\top
\,\mathrm{d}s \rrVert
\stoch0$,
\item[(ii)]
$\sum_{k=1}^{\lfloor nT \rfloor}
\EE ( \|\bU^{(n)}_k\|^2 \bbone_{\{\|\bU^{(n)}_k\| > \theta\}}
|\cF^{(n)}_{k-1}  )
\stoch0$
for all $\theta>0$,
\end{enumerate}
where $\stoch$ denotes convergence in probability.
Then $\bcU^{(n)} \distr\bcU$, as $n\to\infty$.
\end{Thm}

Note that in (i) of Theorem~\ref{Conv2DiffThm}, $\|\cdot\|$ denotes
a matrix norm, while in (ii) it denotes a vector norm.

\section*{Acknowledgements}

The authors have been partially supported by the Hungarian Chinese
Intergovernmental S\&T Cooperation Programme for 2011--2013 under Grant
No. 10-1-2011-0079.
M. Isp\'any has been partially supported by the
T\'AMOP-4.2.2.C-11/1/KONV-2012-0001 project.
The project has been supported by the European Union, co-financed by the
European Social Fund.
K. K\"ormendi was supported by the European Union and the State of Hungary,
co-financed by the European Social Fund in the framework of
T\'AMOP 4.2.4. A/2-11-1-2012-0001 ``National Excellence Program''.
G. Pap has been partially supported by the Hungarian Scientific
Research Fund
under Grant No.~OTKA T-079128.

\end{appendix}


%suskaldyti doi

% imsref loaded by jurgita.kaciuliene, 2014-01-10 11:41:35

\printhistory


\begin{thebibliography}{22}
% Style name=bej, version=1.0, label_style=nolabel, sorting_style=complex, cfg=None, language=None.


%b1 ###
\bibitem{AN}
\begin{bbook}[mr]
\bauthor{\bsnm{Athreya},~\bfnm{Krishna~B.}\binits{K.B.}} \AND
\bauthor{\bsnm{Ney},~\bfnm{Peter~E.}\binits{P.E.}}
(\byear{1972}).
\btitle{Branching Processes}.
\blocation{New York}:
\bpublisher{Springer}.
\bid{mr={0373040}}
\end{bbook}
\bptok{imsref}%
% NOT OUTPUTED:
%   fpage = xi+287
\endbibitem

%b2 ###
\bibitem{BarIspPap2}
\begin{bmisc}[auto:STB|2014/01/06|10:16:28]
\bauthor{\bsnm{Barczy},~\bfnm{M.}\binits{M.}},
\bauthor{\bsnm{Isp{\'a}ny},~\bfnm{M.}\binits{M.}} \AND
\bauthor{\bsnm{Pap},~\bfnm{G.}\binits{G.}}
\bhowpublished{(2012). Asymptotic behavior of CLS estimators
for unstable INAR(2) models. Available at \url{http://arxiv.org/abs/1202.1617}}.
\end{bmisc}
\bptok{imsref}%
% NOT OUTPUTED:
%   sortkey = Barczy(2012
\endbibitem

%b3 ###
\bibitem{BarIspPap0}
\begin{barticle}[mr]
\bauthor{\bsnm{Barczy},~\bfnm{M.}\binits{M.}},
\bauthor{\bsnm{Isp{\'a}ny},~\bfnm{M.}\binits{M.}} \AND
\bauthor{\bsnm{Pap},~\bfnm{G.}\binits{G.}}
(\byear{2011}).
\btitle{Asymptotic behavior of unstable {${\rm INAR}(p)$} processes}.
\bjournal{Stochastic Process. Appl.}
\bvolume{121}
\bpages{583--608}.
\bid{doi={10.1016/j.spa.2010.11.005}, issn={0304-4149}, mr={2763097}}
\end{barticle}
\bptok{imsref}%
% NOT OUTPUTED:
%   issn = 0304-4149
%   url = http://dx.doi.org/10.1016/j.spa.2010.11.005
%   number = 3
%   coden = STOPB7
%   fjournal = Stochastic Processes and their Applications
\endbibitem

%b4 ###
\bibitem{Gut}
\begin{bbook}[mr]
\bauthor{\bsnm{Guttorp},~\bfnm{Peter}\binits{P.}}
(\byear{1991}).
\btitle{Statistical Inference for Branching Processes}.
\bseries{Wiley Series in Probability and Mathematical Statistics}.
\blocation{New York}:
\bpublisher{Wiley}.
\bid{mr={1254434}}
\end{bbook}
\bptok{imsref}%
% NOT OUTPUTED:
%   isbn = 0-471-82291-4
%   fpage = xii+211
\endbibitem

%b5 ###
\bibitem{HalYao}
\begin{barticle}[mr]
\bauthor{\bsnm{Hall},~\bfnm{Peter}\binits{P.}} \AND
\bauthor{\bsnm{Yao},~\bfnm{Qiwei}\binits{Q.}}
(\byear{2003}).
\btitle{Inference in {ARCH} and {GARCH} models with heavy-tailed errors}.
\bjournal{Econometrica}
\bvolume{71}
\bpages{285--317}.
\bid{doi={10.1111/1468-0262.00396}, issn={0012-9682}, mr={1956860}}
\end{barticle}
\bptok{imsref}%
% NOT OUTPUTED:
%   issn = 0012-9682
%   url = http://dx.doi.org/10.1111/1468-0262.00396
%   number = 1
%   coden = ECMTA7
%   fjournal = Econometrica. Journal of the Econometric Society
\endbibitem

%b6 ###
\bibitem{Ham}
\begin{bbook}[mr]
\bauthor{\bsnm{Hamilton},~\bfnm{James~D.}\binits{J.D.}}
(\byear{1994}).
\btitle{Time Series Analysis}.
\blocation{Princeton, NJ}:
\bpublisher{Princeton Univ. Press}.
\bid{mr={1278033}}
\end{bbook}
\bptok{imsref}%
% NOT OUTPUTED:
%   isbn = 0-691-04289-6
%   fpage = xvi+799
\endbibitem

%b7 ###
\bibitem{HJ}
\begin{bbook}[mr]
\bauthor{\bsnm{Horn},~\bfnm{Roger~A.}\binits{R.A.}} \AND
\bauthor{\bsnm{Johnson},~\bfnm{Charles~R.}\binits{C.R.}}
(\byear{1985}).
\btitle{Matrix Analysis}.
\blocation{Cambridge}:
\bpublisher{Cambridge Univ. Press}.
\bid{mr={0832183}}
\end{bbook}
\bptok{imsref}%
% NOT OUTPUTED:
%   isbn = 0-521-30586-1
%   fpage = xiii+561
\endbibitem

%b8 ###
\bibitem{IKP}
\begin{bmisc}[auto:STB|2014/01/06|10:16:28]
\bauthor{\bsnm{Isp{\'a}ny},~\bfnm{M.}\binits{M.}},
\bauthor{\bsnm{K{\"o}rmendi},~\bfnm{K.}\binits{K.}} \AND
\bauthor{\bsnm{Pap},~\bfnm{G.}\binits{G.}}
\bhowpublished{(2012). Asymptotic behavior of CLS estimators for $2$-type critical Galton--Watson
processes with immigration. Available at \url{http://arxiv.org/abs/1210.8315}}.
\end{bmisc}
\bptok{imsref}%
% NOT OUTPUTED:
%   sortkey = Ispany(2012
\endbibitem

%b9 ###
\bibitem{IspPap2}
\begin{bmisc}[auto:STB|2014/01/06|10:16:28]
\bauthor{\bsnm{Isp{\'a}ny},~\bfnm{M.}\binits{M.}} \AND
\bauthor{\bsnm{Pap},~\bfnm{G.}\binits{G.}}
\bhowpublished{(2012). Asymptotic behavior of critical primitive multi-type branching processes with
immigration. Available at \url{http://arxiv.org/abs/1205.0388}}.
\end{bmisc}
\bptok{imsref}%
% NOT OUTPUTED:
%   sortkey = Ispany(2012
\endbibitem

%b10 ###
\bibitem{IspPap}
\begin{barticle}[mr]
\bauthor{\bsnm{Isp{\'a}ny},~\bfnm{M.}\binits{M.}} \AND
\bauthor{\bsnm{Pap},~\bfnm{G.}\binits{G.}}
(\byear{2010}).
\btitle{A note on weak convergence of random step processes}.
\bjournal{Acta Math. Hungar.}
\bvolume{126}
\bpages{381--395}.
\bid{doi={10.1007/s10474-009-9099-5}, issn={0236-5294}, mr={2629664}}
\end{barticle}
\bptok{imsref}%
% NOT OUTPUTED:
%   issn = 0236-5294
%   url = http://dx.doi.org/10.1007/s10474-009-9099-5
%   number = 4
%   fjournal = Acta Mathematica Hungarica
\endbibitem

%b11 ###
\bibitem{JSh}
\begin{bbook}[mr]
\bauthor{\bsnm{Jacod},~\bfnm{Jean}\binits{J.}} \AND
\bauthor{\bsnm{Shiryaev},~\bfnm{Albert~N.}\binits{A.N.}}
(\byear{2003}).
\btitle{Limit Theorems for Stochastic Processes},
\bedition{2nd} ed.
\bseries{Grundlehren der Mathematischen Wissenschaften [Fundamental Principles of Mathematical Sciences]}
\bvolume{288}.
\blocation{Berlin}:
\bpublisher{Springer}.
\bid{mr={1943877}}
\end{bbook}
\bptok{imsref}%
% NOT OUTPUTED:
%   isbn = 3-540-43932-3
%   fpage = xx+661
\endbibitem

%b12 ###
\bibitem{K}
\begin{bbook}[mr]
\bauthor{\bsnm{Kallenberg},~\bfnm{Olav}\binits{O.}}
(\byear{1997}).
\btitle{Foundations of Modern Probability}.
\bseries{Probability and Its Applications (New York)}.
\blocation{New York}:
\bpublisher{Springer}.
\bid{mr={1464694}}
\end{bbook}
\bptok{imsref}%
% NOT OUTPUTED:
%   isbn = 0-387-94957-7
%   fpage = xii+523
\endbibitem

%b13 ###
\bibitem{KesSti1}
\begin{barticle}[mr]
\bauthor{\bsnm{Kesten},~\bfnm{H.}\binits{H.}} \AND
\bauthor{\bsnm{Stigum},~\bfnm{B.~P.}\binits{B.P.}}
(\byear{1966}).
\btitle{A limit theorem for multidimensional {G}alton--{W}atson processes}.
\bjournal{Ann. Math. Statist.}
\bvolume{37}
\bpages{1211--1223}.
\bid{issn={0003-4851}, mr={0198552}}
\end{barticle}
\bptok{imsref}%
% NOT OUTPUTED:
%   issn = 0003-4851
%   fjournal = Annals of Mathematical Statistics
\endbibitem

%b14 ###
\bibitem{MikStr}
\begin{barticle}[mr]
\bauthor{\bsnm{Mikosch},~\bfnm{Thomas}\binits{T.}} \AND
\bauthor{\bsnm{Straumann},~\bfnm{Daniel}\binits{D.}}
(\byear{2002}).
\btitle{Whittle estimation in a heavy-tailed {$\rm GARCH(1,1)$} model}.
\bjournal{Stochastic Process. Appl.}
\bvolume{100}
\bpages{187--222}.
\bid{doi={10.1016/S0304-4149(02)00097-2}, issn={0304-4149}, mr={1919613}}
\end{barticle}
\bptok{imsref}%
% NOT OUTPUTED:
%   issn = 0304-4149
%   url = http://dx.doi.org/10.1016/S0304-4149(02)00097-2
%   coden = STOPB7
%   fjournal = Stochastic Processes and their Applications
\endbibitem

%b15 ###
\bibitem{MR}
\begin{bbook}[mr]
\bauthor{\bsnm{Musiela},~\bfnm{Marek}\binits{M.}} \AND
\bauthor{\bsnm{Rutkowski},~\bfnm{Marek}\binits{M.}}
(\byear{1997}).
\btitle{Martingale Methods in Financial Modelling}.
\bseries{Applications of Mathematics (New York)}
\bvolume{36}.
\blocation{Berlin}:
\bpublisher{Springer}.
\bid{mr={1474500}}
\end{bbook}
\bptok{imsref}%
% NOT OUTPUTED:
%   isbn = 3-540-61477-X
%   fpage = xii+512
\endbibitem

%b16 ###
\bibitem{Q}
\begin{barticle}[mr]
\bauthor{\bsnm{Quine},~\bfnm{M.~P.}\binits{M.P.}}
(\byear{1970}).
\btitle{The multi-type {G}alton--{W}atson process with immigration}.
\bjournal{J. Appl. Probability}
\bvolume{7}
\bpages{411--422}.
\bid{issn={0021-9002}, mr={0263168}}
\end{barticle}
\bptok{imsref}%
% NOT OUTPUTED:
%   issn = 0021-9002
%   fjournal = Journal of Applied Probability
\endbibitem

%b17 ###
\bibitem{RevYor}
\begin{bbook}[auto:STB|2014/01/06|10:16:28]
\bauthor{\bsnm{Revuz},~\bfnm{D.}\binits{D.}} \AND
\bauthor{\bsnm{Yor},~\bfnm{M.}\binits{M.}}
(\byear{2001}).
\btitle{Continuous Martingales and Brownian Motion},
\bedition{3rd} ed.,
\bnote{corrected 2nd printing}.
\blocation{Berlin}:
\bpublisher{Springer}.
\bnote{\MR{1725357}}
\end{bbook}
\bptok{imsref}%
\endbibitem

%b18 ###
\bibitem{SheSri}
\begin{barticle}[mr]
\bauthor{\bsnm{Shete},~\bfnm{Sanjay}\binits{S.}} \AND
\bauthor{\bsnm{Sriram},~\bfnm{T.~N.}\binits{T.N.}}
(\byear{2003}).
\btitle{A note on estimation in multitype supercritical branching processes with immigration}.
\bjournal{Sankhy\=a}
\bvolume{65}
\bpages{107--121}.
\bid{issn={0972-7671}, mr={2016780}}
\end{barticle}
\bptok{imsref}%
% NOT OUTPUTED:
%   issn = 0972-7671
%   number = 1
%   fjournal = Sankhy\=a. The Indian Journal of Statistics
\endbibitem

%b19 ###
\bibitem{Tan}
\begin{bbook}[mr]
\bauthor{\bsnm{Tanaka},~\bfnm{Katsuto}\binits{K.}}
(\byear{1996}).
\btitle{Time Series Analysis: Nonstationary and Noninvertible Distribution Theory}.
\bseries{Wiley Series in Probability and Statistics}.
\blocation{New York}:
\bpublisher{Wiley}.
\bid{mr={1397269}}
\end{bbook}
\bptok{imsref}%
% NOT OUTPUTED:
%   isbn = 0-471-14191-7
%   fpage = x+623
\endbibitem

%b20 ###
\bibitem{WW}
\begin{barticle}[mr]
\bauthor{\bsnm{Wei},~\bfnm{C.~Z.}\binits{C.Z.}} \AND
\bauthor{\bsnm{Winnicki},~\bfnm{J.}\binits{J.}}
(\byear{1989}).
\btitle{Some asymptotic results for the branching process with immigration}.
\bjournal{Stochastic Process. Appl.}
\bvolume{31}
\bpages{261--282}.
\bid{doi={10.1016/0304-4149(89)90092-6}, issn={0304-4149}, mr={0998117}}
\end{barticle}
\bptok{imsref}%
% NOT OUTPUTED:
%   issn = 0304-4149
%   url = http://dx.doi.org/10.1016/0304-4149(89)90092-6
%   number = 2
%   coden = STOPB7
%   fjournal = Stochastic Processes and their Applications
\endbibitem

%b21 ###
\bibitem{WW2}
\begin{barticle}[mr]
\bauthor{\bsnm{Wei},~\bfnm{C.~Z.}\binits{C.Z.}} \AND
\bauthor{\bsnm{Winnicki},~\bfnm{J.}\binits{J.}}
(\byear{1990}).
\btitle{Estimation of the means in the branching process with immigration}.
\bjournal{Ann. Statist.}
\bvolume{18}
\bpages{1757--1773}.
\bid{doi={10.1214/aos/1176347876}, issn={0090-5364}, mr={1074433}}
\end{barticle}
\bptok{imsref}%
% NOT OUTPUTED:
%   issn = 0090-5364
%   url = http://dx.doi.org/10.1214/aos/1176347876
%   number = 4
%   coden = ASTSC7
%   fjournal = The Annals of Statistics
\endbibitem

%b22 ###
\bibitem{Win}
\begin{barticle}[mr]
\bauthor{\bsnm{Winnicki},~\bfnm{J.}\binits{J.}}
(\byear{1991}).
\btitle{Estimation of the variances in the branching process with immigration}.
\bjournal{Probab. Theory Related Fields}
\bvolume{88}
\bpages{77--106}.
\bid{doi={10.1007/BF01193583}, issn={0178-8051}, mr={1094078}}
\end{barticle}
\bptok{imsref}%
% NOT OUTPUTED:
%   issn = 0178-8051
%   url = http://dx.doi.org/10.1007/BF01193583
%   number = 1
%   coden = PTRFEU
%   fjournal = Probability Theory and Related Fields
\endbibitem

\end{thebibliography}
\end{document}